\documentclass[12pt,leqno,twoside]{article} 
\usepackage{bezier,eufrak}
\input{rogg_kue.sty}
\newcommand{\enger}{\setlength{\arraycolsep}{1pt}
                           \renewcommand{\arraystretch}{0.5} }
\newcommand{\weiter}{\setlength{\arraycolsep}{3pt}
                           \renewcommand{\arraystretch}{1.2} }
\weiter

\begin{document}
\title{On the cyclotomic Dedekind embedding and the cyclic Wedderburn embedding}
\author{M. K\"unzer, H.\ Weber}
\date{March 20, 2001}
\maketitle

\begin{small}
\begin{quote}
\begin{center}{\bf Abstract}\end{center}\vspace*{2mm}
Let $n\geq 1$ and let $p$ be a prime. Expand $j\in [0,p^n-1]\ohne (p)$ $p$-adically as $j = \sum_{s\geq 0} a_s p^s$ with $a_s\in [0,p-1]$. The
$\#([0,j]\ohne (p))$th $\Z_{(p)}[\zeta_{p^n}]$-linear elementary divisor of the cyclotomic Dedekind embedding
\[
\Z_{(p)}[\zeta_{p^n}]\ts_{\sZ_{(p)}} \Z_{(p)}[\zeta_{p^n}] \hra \prod_{i\in (\sZ/p^n)^\ast} \Z_{(p)}[\zeta_{p^n}]
\]
has valuation 
\[
-1 + \sumd{s\geq 0} \left(a_s(s+1) - a_{s+1}(s+2)\right) p^s
\]
at $1-\zeta_{p^n}$. There is a similar result for the related cyclic Wedderburn embedding.
\end{quote}
\end{small}

\setcounter{section}{-1}

\section{Introduction}

\setcounter{subsection}{-1}

\subsection{Contents}

\vspace*{-10mm}

\begin{footnotesize}
\parskip0.4ex
\tableofcontents
\parskip1.2ex
\end{footnotesize}

\subsection{Objectives}
In this note, we consider an extension of discrete valuation rings $S\tm T$, with $s\in S$ and $t\in T$ generating the respective maximal ideals. Let $K$ be the field of fractions of $S$, 
let $L$ be the field of fractions of $T$ and assume $L$ over $K$ to be finite galois of Galois group $G$ of order $m = |G|$. Moreover, assume $T = S[\theta]$ for some $\theta\in T$. We restrict the 
Dedekind isomorphism
\[
\begin{array}{rcl}
L\ts_K L & \lraiso  & \prodd{\sigma\in G} L \\
x\ts y   & \lramaps & (xy^\sigma)_{\sigma\in G}, \\
\end{array}
\]
the injectivity of which ensues from Dedekind's Lemma, to the locally integral situation
\[
T\ts_S T\;\hra\; \prodd{\sigma\in G} T.
\]
This embedding is no longer an isomorphism in general; in fact, the square of its $T$-linear determinant is just the discriminant of $T$ over $S$. Being a $T$-linear map, we can ask for its 
elementary divisors, or, more precisely, for an elementary divisor diagonalization. 

Once such a diagonalization obtained, we may use the image of $T\ts_S T$ inside $\prod_{\sigma\in G} T$, instead of $T\ts_S T$ itself, to calculate Hochschild (co)homology of $T$. For instance, we 
may ask for these (co)homology groups with coefficients in $T$, equipped with a galois twisted $T$-$T$-bimodule structure (the untwisted case being due to {\sc M.\ Larsen} and 
{\sc A.\ Lindenstrauss} [LL 92, 1.6.2] and to the Buenos Aires Cyclic Homology group [BACH 91, prop.\ 1.3]).

The cyclotomic Dedekind embedding shall be studied in more detail, where $T = \Z_{(p)}[\zeta_{p^n}]$ and $S = \Z_{(p)}$, for $n\geq 1$, $p$ a prime, and where, in general for $m\geq 1$,
$\zeta_m$ denotes a primitive $m$th root of unity over $\Q$.

This cyclotomic case leads to a related question, which is to give an elementary divisor diagonalization for the cyclic Wedderburn embedding. Given $m\geq 1$, and letting $c_m$ denote a 
generator of the cyclic group $C_m$ of $m$ elements, this embedding results from restricting the rational Wedderburn isomorphism
\[
\begin{array}{rcl}
\Q(\zeta_m)C_m & \lraiso  & \prodd{j\in [0,m-1]} \Q(\zeta_m) \\
c_m            & \lramaps & (\zeta_m^j)_{j\in [0,m-1]}, \\
\end{array}
\]
to the integral situation
\[
\Z[\zeta_m] C_m \;\hra\; \prodd{j\in [0,m-1]} \Z[\zeta_m],
\]
which is likewise no longer an isomorphism for $m > 1$, the square of its $\Z[\zeta_m]$-linear determinant being $\pm m^m$. Now the cyclotomic Dedekind embedding may be viewed, considered as a 
diagram, as a quotient of the corresponding cyclic Wedderburn embedding.

This consideration in turn gives rise to the question for a closed description of the image of the cyclic Wedderburn embedding in the absolute case, i.e.\ of
\[
\begin{array}{rcl}
\Z C_m & \hra     & \prodd{d\, |\, m} \Z[\zeta_d] \\
c_m    & \lramaps & (\zeta_d)_{d\,|\, m}. \\
\end{array}
\]

\subsection{Results}

\subsubsection{Hochschild (co)homology with twisted coefficients}

We choose a total ordering on $G$, $(\sigma_0,\dots,\sigma_{m-1})$, in such a way that $\sigma_0 = 1_T$, and such that
\[
\sumd{i\in [0,j-1]} v_t(\theta^{\sigma_j} - \theta^{\sigma_i}) \leq \sumd{i\in [0,j-1]} v_t(\theta^{\sigma_k} - \theta^{\sigma_i}) 
\]
for each $j\in [0,m-1]$ and each $k\in [j+1,m-1]$, where $v_t$ denotes the $t$-adic valuation. Let $T_i$ denote $T$, considered as a $T$-$T$-bimodule with untwisted left multiplication and 
$\sigma_i$-twisted right multiplication by $T$. Let $\phi := \sum_{k\in [1,m-1]} v_t(\theta - \theta^{\sigma_k})$.
 
The Hochschild homology of $T$, over the ground ring $S$ and with coefficients in $T_i$, $i\in [0,m-1]$, is given in dimension $j\geq 0$ by
\[
H_j(T,T_0;S) \;\iso\;
\left\{
\begin{array}{ll}
T           & \mb{for $j = 0$} \\
T/t^{\phi}T & \mb{for $j$ odd} \\
0           & \mb{for $j$ even, $j\geq 2$} \\
\end{array}
\right.
\]
and by 
\[
H_j(T,T_i;S) \;\iso\;
\left\{
\begin{array}{ll}
0                                & \mb{for $j$ odd} \\
T/(\theta^{\sigma_i} - \theta) T & \mb{for $j$ even} \\
\end{array}
\right.
\]
for $i\in [1,m-1]$. The Hochschild cohomology is given by
\[
H^j(T,T_0;S) \;\iso\;
\left\{
\begin{array}{ll}
T           & \mb{for $j = 0$} \\
0           & \mb{for $j$ odd} \\
T/t^{\phi}T & \mb{for $j$ even, $j\geq 2$} \\
\end{array}
\right.
\]
and by 
\[
H^j(T,T_i;S) \;\iso\;
\left\{
\begin{array}{ll}
T/(\theta^{\sigma_i} - \theta) T & \mb{for $j$ odd} \\
0                                & \mb{for $j$ even} \\
\end{array}
\right.
\]
for $i\in [1,m-1]$ (\ref{CorPR3}). 

\subsubsection{Elementary divisors of the cyclotomic Dedekind embedding}

Let $n\geq 1$ and let $p$ be a prime. Expand $j\in [0,p^n-1]\ohne (p)$ $p$-adically as $j = \sum_{s\geq 0} a_s p^s$ with $a_s\in [0,p-1]$. The
$\#([0,j]\ohne (p))$th $\Z_{(p)}[\zeta_{p^n}]$-linear elementary divisor of the cyclotomic Dedekind embedding
\[
\begin{array}{rllcl}
\Z_{(p)}[\zeta_{p^n}] & \ts_{\sZ_{(p)}} & \Z_{(p)}[\zeta_{p^n}] & \hrah{\delr_{p^n}} & \prod_{i\in (\sZ/p^n)^\ast} \Z_{(p)}[\zeta_{p^n}] \\
\zeta_{p^n}^k         & \ts             & \zeta_{p^n}^l         & \lramaps           & (\zeta_{p^n}^{k+il}), \\
\end{array}
\]
where $k,l\in [0,p^{n-1}(p-1)-1]$, has valuation 
\[
-1 + \sumd{s\geq 0} \left(a_s(s+1) - a_{s+1}(s+2)\right) p^s
\]
at $1-\zeta_{p^n}$ (\ref{PropVand9}).

\subsubsection{The cyclic Wedderburn embedding}

Expand $j\in [0,p^n-1]$ $p$-adically as $j = \sum_{s\geq 0} a_s p^s$ with $a_s\in [0,p-1]$. The $(j+1)$st $\Z_{(p)}[\zeta_{p^n}]$-linear elementary divisor of the localized cyclic Wedderburn 
embedding
\[
\begin{array}{rcl}
\Z_{(p)}[\zeta_{p^n}] C_{p^n} & \hrah{\omega_{p^n}} & \prod_{j\in\sZ/p^n} \Z_{(p)}[\zeta_{p^n}] \\
c_{p^n}                       & \lramaps            & (\zeta_{p^n}^j)_{j\in\sZ/p^n} \\
\end{array}
\]
has valuation
\[
\sumd{s\geq 0} (a_s - a_{s+1})(s+1)p^s
\]
at $1-\zeta_{p^n}$ (\ref{PropVand5}). We can render this result a bit more precise and a bit more general, in that we can diagonalize the cyclic Wedderburn embedding 
\[
\begin{array}{rcl}
\Z[\zeta_m] C_m & \hrah{\omega_m} & \prod_{j\in\sZ/m} \Z[\zeta_m] \\
c_m             & \lramaps        & (\zeta_m^j)_{j\in\sZ/m} \\
\end{array}
\]
for $m\geq 1$, regardless whether $\Z[\zeta_m]$ is a principal ideal domain (\ref{ThSaal5}). For $l\in [0,m-1]$, the $(l+1)$th diagonal entry becomes
\[
\fracd{m \zeta_m^{(l^2)}}{\prod_{j\in [1,l]} (1-\zeta_m^j)}.
\]

\subsubsection{The absolute cyclic Wedderburn embedding}

It suffices to consider the primary parts of $m$ separately, i.e.\ we may assume $m = p^n$ (cf.\ \ref{PropCP9_2}). In this case, the absolute cyclic Wedderburn embedding is given by
\[
\begin{array}{rcl}
\Z C_{p^n} & \hrah{\omega_{\ssZ,p^n}} & \prodd{i\in [0,n]} \Z[\zeta_{p^i}] \\
c_{p^n}    & \lramaps                 & (\zeta_{p^i})_{i\in [0,n]}. \\
\end{array}
\]
We derive from the pullback of {\sc Kervaire} and {\sc Murthy} [KM 77] the following triangular system of ties, i.e.\ of congruences of tuple entries, that describe its image.

Let $\phi$ denote Euler's function. For $m\geq 1$ and $j\in\Z$, we let $[j]_m\in [0,m-1]$ be such that $[j]_m\con_m j$. Let $\dell$ denote Kronecker's delta. 
The image of the absolute Wedderburn embedding is given by
\[
\begin{array}{rcl}
(\Z C_{p^n})\omega_{\sZ,p^n} 
& =   & \left\{\rule[7mm]{0mm}{0mm}\right.\left(\sumd{j\in [0,\phi(p^i) - 1]} x_{i,j}\zeta_{p^i}^j\right)_{\!\! i\in [0,n]}\hspace*{-7mm},\; x_{i,j}\in\Z\;\left|\rule[7mm]{0mm}{0mm}\right. \\
&     & \mb{for $l\in [1,n]$ and $j\in [0,\phi(p^{n-l})-1]$ we have }\; x_{n-l,j} \con_{p^l} \sumd{i\in [0,l-1]} p^{l-1-i} \cdot\\
&     & \cdot \sumd{k\in [1,p-1]} \left(x_{n-i,j - p^{n-l} + kp^{n-1-i}} - (1-\dell_{l,n})x_{n-i,[j]_{p^{n-l-1}} - p^{n-l-1} + kp^{n-1-i}}\right)
        \left.\rule[7mm]{0mm}{0mm}\right\} \\
& \tm & \prodd{i\in [0,n]}\Z[\zeta_{p^i}]. \\
\end{array}
\]
The elementary divisors of $\omega_{\sZ,p^n}$ over $\Z$ are given by $p^i$ with multiplicity $\phi(p^{n-i})$ for $i\in [0,n]$ (\ref{ThKM5}).

{\sc Kleinert} gives a system of ties that describes the image of the absolute Wedderburn embedding $\Z C_m\lraa{\omega_{\ssZ,m}}\prod_{d|m} \Z[\zeta_d]$ in terms of certain 
prime ideals of the rings $\Z[\zeta_d]$, $d|m$, in case $m$ is a positive squarefree integer [Kl 81, p.\ 550]. In loc.\ cit., this description is used as a tool to 
study units in dihedral group rings.

\subsection{Acknowledgements}

We would like to thank {\sc T.\ Zink} for a hint how to simplify our elementary divisor calculations, see subsections \ref{SubSecDiagVand}, \ref{SubSecElVand}. We would like to thank 
{\sc J.\ Stokman} for the deduction of the diagonalization of the Wedderburn embedding from the $q$-Saalsch\"utz Theorem, see subsection \ref{SubSecSaal}. We would like to thank several colleagues 
for comments on an earlier version. The second author would like to thank {\sc P.\ Littelmann} for kind hospitality in Strasbourg, where this earlier version has been written, and where he has 
been supported by the EU TMR-network `Algebraic Lie Representations', grant no.\ ERB FMRX-CT97-0100.

\vspace*{6mm}

\pagebreak[2]           
\begin{footnotesize}
\begin{Notation}
\label{NotIF0}
\Absit\rm
\begin{itemize}
\item[(o)] Composition of maps is written on the right, $\lraa{a}\lraa{b} = \lraa{ab}$. 
\item[(i)] If $\phi$ is a map of $A$-modules, $A$ a commutative ring, and $\pfk$ a prime ideal of $A$, we sometimes denote the localization $\phi_\pfk$ merely by $\phi$. If $m\geq 1$, we
denote the ring of $m\ti m$-matrices with entries in $A$ by $(A)_m$.
\item[(ii)] If $x$, $y$ are elements of some set, we let $\dell_{x,y} := 1$ if $x = y$ and $\dell_{x,y} := 0$ if $x\neq y$. 
\item[(iii)] For integers $a$, $b$, we denote $[a,b] := \{x\in\Z\; |\; a\leq x\leq b\}$ and $[a,b[\, := \{x\in\Z\; |\; a\leq x < b\}$.
\item[(iv)] For $a, b\geq 0$, we let the binomial coefficient $\smatze{a}{b}$ be equal to zero if $b > a$.
\item[(v)] Let $\ts = \ts_{\sZ}$ (i.e.\ $\ts = \ts_{\sQ}$ over $\Q$). 
\end{itemize}
Let $m\geq 1$.
\begin{itemize}
\item[(vi)]  For a prime $p$, we let $m[p] := p^{v_p(m)}$ denote the $p$-part of $m$. 
\item[(vii)] Let $\phi$ denote Euler's function, $\phi(m) := m\prod_{p|m} (1 - p^{-1})$. Let $\Phi_m(X)$ denote
the $m$-th cyclotomic polynomial, i.e.\ the irreducible factor of $X^m - 1\in\Q[X]$ that does not divide $X^{m'}-1$ for any $m'|m$, $m'\neq m$. Hence $\deg\Phi_m(X) = \phi(m)$.
\item[(viii)] Let $\zeta_m$ denote a primitive $m$th root of unity over $\Q$, with minimal polynomial $\mu_{\zeta_m,\sQ}(T) = \Phi_m(T)$. Let $c_m$ denote a generator of the cyclic group $C_m$ of 
$m$ elements.
\end{itemize}
\end{Notation}
\end{footnotesize}

\section{Vandermonde}
\label{SecVand}

\begin{quote}
\begin{footnotesize}
We diagonalize the polynomial Vandermonde matrix rationally. Specializing under the assumption that the diagonalizing matrices become integral, we obtain an elementary divisor form, i.e.\ a 
diagonal matrix in which each subsequent diagonal entry is divisible by its predecessor. 
\end{footnotesize}
\end{quote}

\subsection{Diagonalization of Vandermonde matrices}
\label{SubSecDiagVand}

\begin{quote}
\begin{footnotesize}
A hint of {\sc T.\ Zink} how to simplify our previous efforts led to the approach in this and in the next subsection. 
\end{footnotesize}
\end{quote}

\begin{Notation}
\label{NotDiag1}\rm
Let $m\geq 1$. Let $x = (x_0,\dots, x_{m-1})$ be a tuple of indeterminates, let $a,b\in [0,m-1]$ and let $d\in\Z$. We consider in $\Z[x_0,\dots,x_{m-1}]$ the symmetric polynomials
\[
\begin{array}{rcl}
P_{d,[a,b]} & := & 
\left\{
\begin{array}{ll}
\sumd{s_i\geq 0,\; \sum_{i\in [a,b]} s_i \; =\; d} x_a^{s_a} x_{a+1}^{s_{a+1}}\cdots x_b^{s_b} & \mb{if $d\geq 0$}\\
0                                                                                              & \mb{if $d < 0$} \\
\end{array}
\right.\vspace*{4mm} \\
E_{d,[a,b[} & := & 
\left\{
\begin{array}{ll}
\sumd{s_i\in\{0,1\},\; \sum_{i\in [a,b[} s_i \; =\; d} x_a^{s_a} x_{a+1}^{s_{a+1}}\cdots x_{b-1}^{s_{b-1}} & \mb{if $d\geq 0$} \\
0                                                                                                          & \mb{if $d < 0$.} \\
\end{array}
\right. \\
\end{array}
\]
In particular, we have
\[
\begin{array}{rcl}
P_{d,\leer} & = & \dell_{d,0} \\
E_{d,\leer} & = & \dell_{d,0}.\\
\end{array}
\]
Moreover, for $i\in [0,m-1]$ we denote
\[
y_i \; := \; \prod_{k\in [0,i[} (x_i - x_k) = \sum_{j\in [0,i]} (-1)^{i-j} E_{i-j,[0,i[} x_i^j,
\]
considered as element of $\Z[x_0,\dots, x_{m-1}]$. For $i,j\in [0,m-1]$, we denote the Lagrange interpolation function by
\[
L_{i,j} \; :=\; 
\left\{
\begin{array}{ll}
\fracd{\prod_{k\in [0,j[\ohne\{ i\}} (x_j - x_k)}{\prod_{k\in [0,j[\ohne\{ i\}} (x_i - x_k)} & \mb{if $i < j$} \\
0                                                                                            & \mb{if $i\geq j$} \\
\end{array}
\right.,
\]
and as a variant
\[
M_{i,j} := 
\fracd{\prod_{k\in [0,i[} (x_j - x_k)}{\prod_{k\in [0,i[} (x_i - x_k)},
\]
considered as elements of $\Q(x_0,\dots,x_{m-1})$. These polynomials furnish the matrices
\[
\begin{array}{rcl}
I & := & (\dell_{i,j})_{i\ti j\in [0,m-1]\ti [0,m-1]}\\
V_x & := & (x_j^i)_{i\ti j\in [0,m-1]\ti [0,m-1]} \\
L_x & := & (L_{i,j})_{i\ti j\in [0,m-1]\ti [0,m-1]} \\
M_x & := & (M_{i,j})_{i\ti j\in [0,m-1]\ti [0,m-1]} \\
P_x & := & (P_{i-j,[0,j]})_{i\ti j\in [0,m-1]\ti [0,m-1]} \\
E_x & := & ((-1)^{i-j} E_{i-j,[0,i[})_{i\ti j\in [0,m-1]\ti [0,m-1]} \\
Y_x & := & (\dell_{i,j} y_i)_{i\ti j\in [0,m-1]\ti [0,m-1]} \\
\end{array}
\]
in $(\Q(x_0,\dots,x_{m-1}))_m$.
\end{Notation}

\begin{Lemma}
\label{LemDiag2}
For $0\leq a\leq b\leq m$ and $d\in\Z$, we have
\[
\begin{array}{rclcl}
P_{d,[a,b]} & = & P_{d,[a+1,b]} + x_a P_{d-1,[a,b]}     & = & P_{d,[a,b-1]} + x_b P_{d-1,[a,b]} \\
E_{d,[a,b[} & = & E_{d,[a+1,b[} + x_a E_{d-1,[a+1,b[} & = & E_{d,[a,b-1[} + x_b E_{d-1,[a,b-1[}, \\
\end{array}
\]
where for the second equation, we stipulate in addition that $a < b$.
\end{Lemma}

\begin{Lemma}
\label{LemDiag3}
For $i,k\in\Z$ and $0\leq a\leq b < c\leq m$, we obtain
\[
\sum_{j\in \sZ} (-1)^j E_{i-j,[a,c[} P_{j-k,[a,b]} = (-1)^k E_{i-k,[b+1,c[}.
\]

\rm
In fact,
\[
\begin{array}{rcl}
\sum_{j\in \sZ} (-1)^j E_{k-j,[a,c[} P_{j-i,[a,b]} 
& \aufgl{(\ref{LemDiag2})} & \sum_{j\in \sZ} (-1)^j E_{k-j,[a+1,c[} P_{j-i,[a,b]} \\ 
& + & \sum_{j\in \sZ} (-1)^j x_a E_{k-j-1,[a+1,c[} P_{j-i,[a,b]} \\ 
& = & \sum_{j\in \sZ} (-1)^j E_{k-j,[a+1,c[} P_{j-i,[a,b]} \\ 
& + & \sum_{j\in \sZ} (-1)^{j-1} E_{k-j,[a+1,c[} x_a P_{j-1-i,[a,b]} \\ 
& \aufgl{(\ref{LemDiag2})} & \sum_{j\in \sZ} (-1)^j E_{k-j,[a+1,c[} P_{j-i,[a,b]} \\ 
& + & \sum_{j\in \sZ} (-1)^{j-1} E_{k-j,[a+1,c[} P_{j-i,[a,b]} \\ 
& - & \sum_{j\in \sZ} (-1)^{j-1} E_{k-j,[a+1,c[} P_{j-i,[a+1,b]} \\ 
& = & \sum_{j\in \sZ} (-1)^j E_{k-j,[a+1,c[} P_{j-i,[a+1,b]}. \\
\end{array}
\]
\end{Lemma}

\begin{Lemma}
\label{LemDiag5}
We have $E_x P_x = I$.

\rm
Given $i,k\in [0,m-1]$, we may assume $i > k$ and obtain
\[
\begin{array}{rcl}
(E_x P_x)_{i,k}
& = & \sum_{j\in [0,m-1]} (-1)^{i-j} E_{i-j,[0,i[} P_{j-k,[0,k]} \\
& \aufgl{(\ref{LemDiag3})} & (-1)^{i-k} E_{i-k,[k+1,i[} \\
& = & 0. \\
\end{array}
\]
\end{Lemma}

\begin{Lemma}
\label{LemDiag5_5}
We have $M_x(I-L_x) = I$. 

\rm
Given $i,k\in [0,m-1]$, we may assume $i < k$ and need to show that
\[
M_{i,k} - \sum_{j\in [i,k[} M_{i,j} L_{j,k} = 0.
\]
This expression being an element of $\Q(x_0,\dots,x_{k-1},x_{k+1},x_{m-1})[x_k]$ of degree in $x_k$ less than $k$, we are reduced to plug in $x_k = x_l$ for $l\in [0,k-1]$.
\end{Lemma}

\begin{Lemma}
\label{LemDiag6}
We have $E_x V_x = Y_x M_x$.

\rm
Given $i,k\in [0,m-1]$, we have
\[
\sum_{j\in [0,m-1]} (-1)^{i-j} E_{i-j,[0,i[} x_k^j = \prod_{j\in [0,i[} (x_k - x_j) = y_i M_{i,k}.
\]
\end{Lemma}

\begin{Proposition}
\label{PropDiag7}
We have $E_x V_x(I-L_x) = Y_x$.

\rm
This follows from (\ref{LemDiag5_5}, \ref{LemDiag6}).
\end{Proposition}

\begin{Remark}
\label{RemDiag8}
We have $V_x = P_x Y_x M_x$.

\rm
This follows from (\ref{LemDiag5}, \ref{LemDiag6}).
\end{Remark}

\begin{quote}
\begin{footnotesize}
\begin{Example}
\label{ExDiag9}\rm
Letting $m = 4$, we obtain
\[
\setlength{\arraycolsep}{2pt}
\begin{array}{l}
E_x V_x(I-L_x) \vspace*{2mm}\\
= \left[
\begin{array}{cccc}
\scm 1            & \scm 0                           & \scm 0                   & \scm 0 \\
\scm -x_0         & \scm 1                           & \scm 0                   & \scm 0 \\
\scm x_0 x_1      & \scm -(x_0 + x_1)                & \scm 1                   & \scm 0 \\
\scm -x_0 x_1 x_2 & \scm x_0 x_1 + x_0 x_2 + x_1 x_2 & \scm - (x_0 + x_1 + x_2) & \scm 1 \\
\end{array}
\right] 
\cdot
\left[
\begin{array}{cccc}
\scm 1     & \scm 1     & \scm 1     & \scm 1     \\
\scm x_0   & \scm x_1   & \scm x_2   & \scm x_3   \\
\scm x_0^2 & \scm x_1^2 & \scm x_2^2 & \scm x_3^2 \\
\scm x_0^3 & \scm x_1^3 & \scm x_2^3 & \scm x_3^3 \\
\end{array}
\right] 
\cdot
\left[
\begin{array}{cccc}
\scm 1 & \scm -1 & \scm -\frac{x_2-x_1}{x_0-x_1} & \scm -\frac{(x_3-x_1)(x_3-x_2)}{(x_0-x_1)(x_0-x_2)} \\
\scm 0 & \scm  1 & \scm -\frac{x_2-x_0}{x_1-x_0} & \scm -\frac{(x_3-x_0)(x_3-x_2)}{(x_1-x_0)(x_1-x_2)} \\
\scm 0 & \scm  0 & \scm 1                        & \scm -\frac{(x_3-x_0)(x_3-x_1)}{(x_2-x_0)(x_2-x_1)} \\
\scm 0 & \scm  0 & \scm 0                        & \scm 1                                              \\
\end{array}
\right] \vspace*{2mm}\\
=
\left[
\begin{array}{cccc}
\scm 1 & \scm 0         & \scm 0                  & \scm 0 \\ 
\scm 0 & \scm (x_1-x_0) & \scm 0                  & \scm 0 \\ 
\scm 0 & \scm 0         & \scm (x_2-x_0)(x_2-x_1) & \scm 0 \\ 
\scm 0 & \scm 0         & \scm 0                  & \scm (x_3-x_0)(x_3 - x_1)(x_3 - x_2) \\ 
\end{array}
\right] \;\; = Y_x. \\
\end{array}
\setlength{\arraycolsep}{3pt}
\]
\end{Example}
\end{footnotesize}
\end{quote}

\subsection{Elementary divisors of Vandermonde matrices over discrete valuation rings}
\label{SubSecElVand}

\begin{Setup}\rm
\label{SetupEld0}
Let $T$ be a discrete valuation ring with fraction field $L = \fracfield T$. Let $t\in T$ be a generator of the maximal ideal, and let $v_t$ be the according valuation.
\end{Setup}

\begin{Definition}
\label{DefEld1}\rm
Let $m\geq 1$. Given a tuple $\xi = (\xi_0,\dots,\xi_{m-1})$ of pairwise distinct elements of $T$, we say that $\xi$ is {\it minimally ordered,} if
\[
\sumd{i\in [0,j-1]} v_t(\xi_j - \xi_i) \leq \sumd{i\in [0,j-1]} v_t(\xi_k - \xi_i)
\]
for each $j\in [0,m-1]$ and each $k\in [j+1,m-1]$. Note that any tuple of pairwise distinct elements of $T$ can be reordered (non-uniquely, in general) to a minimally ordered tuple.
\end{Definition}

\begin{Lemma}
\label{LemEld2}
If $\xi$ is minimally ordered, then $L_\xi\in (T)_m$ (cf.\ \ref{NotDiag1}). In particular, the $(i+1)$st $T$-linear elementary divisor of $V_\xi$, where $i\in [0,m-1]$, has valuation
\[
\sum_{j\in [0,i-1]} v_t(\xi_i - \xi_j).
\]

\rm
By (\ref{LemDiag5_5}), $L_\xi$ is contained in $(T)_m$ if and only if $M_\xi$ is contained in $(T)_m$. The assertion on the elementary divisors now ensues from (\ref{PropDiag7}).
\end{Lemma}

\begin{Setup}\rm
\label{SetupEld3}
Let $S\tm T$ be a finite extension of discrete valuation rings, with fraction fields $K := \fracfield S$ and $L := \fracfield T$. Let $s\in S$ and $t\in T$ be generators of the respective maximal ideals, and 
$v_s$ and $v_t$ the respective valuations. Assume $L$ over $K$ to be galois of degree $m = [L:K]$ with Galois group $G$. Assume that 
\[
T = S[\theta] 
\]
for some $\theta\in T$. Choose an ordering 
\[
\begin{array}{rcl}
[0,m-1] & \lraiso  & G \\
i       & \lramaps & \sigma_i \\
\end{array}
\]
such that the tuple $\tau := (\theta^{\sigma_1},\dots,\theta^{\sigma_{m-1}})$ is minimally ordered, starting with $\sigma_0 = 1_T$.
\end{Setup}

\begin{Proposition}
\label{ThEld4}
For $i,j\in [0,m-1]$, we let $L_{i,j}(\tau)$ denote the specialization of $L_{i,j}$ along $x\lramaps\tau$ (cf.\ \ref{NotDiag1}). The following assertions hold.

\begin{itemize}
\item[(i)] The $(i+1)$st $T$-linear elementary divisor, $i\in [0,m-1]$, of the Dedekind embedding 
\[
\begin{array}{rllcl}
T\; & \ts_S & T & \hrah{\delr_{T/S}} & \prod_{\sigma\in G} T \\
x\; & \ts   & y & \lramaps           & (xy^\sigma)_{\sigma\in G} \\
\end{array}
\]
has valuation
\[
\phi_i := \sum_{j\in [0,i-1]} v_t(\theta^{\sigma_i} - \theta^{\sigma_j}).
\]
In particular, $(\phi_i)_{i\in [0,m-1]}$ does neither depend on the choice of $\theta$ nor on the minimal ordering chosen on $\tau$.
\item[(ii)] The image of $T\ts_S T$ under $\delr_{T/S}$, which is an isomorphic copy of $T\ts_S T$, allows a description via ties, i.e.\ via congruences of tuple entries, as
\[
\fbox{$\;\;
\begin{array}{rl}
    & \hspace*{-8mm} (T\ts_S T)\delr_{T/S} \\
=   & \left\{ (\eta_j)_{j\in [0,m-1]}\;\left|\rule[7mm]{0mm}{0mm}\right.\; 
\eta_i - \sumd{j\in [0,i-1]} \eta_j L_{j,i}(\tau) \in T t^{\phi_i}\;\mb{\rm\ for $i\in [0,m-1]$} \right\} \\
\tm & \prod_{j\in [0,m-1]} T. \\
\end{array}
\;\;$}
\]

\item[(iii)] A $T$-linear basis of $(T\ts_S T)\delr_{T/S}$ of triangular shape is given by the following tuple of elements of $\prod_{j\in [0,m-1]} T$.
\[
\left(\left(\prod_{k\in [0,i-1]} (\theta^{\sigma_j} - \theta^{\sigma_k})\right)_{j\in [0,m-1]} \right)_{i\in [0,m-1]}.
\]
\end{itemize}

\rm 
Ad (i). This follows from (\ref{LemEld2}).

Ad (ii). An element $\eta$ of $\prod_{i\in [0,m-1]} T$, considered as a row vector, is contained in the image of $\delr_{T/S}$ if and only if its rational inverse image, written as a row vector
in the $L$-linear basis $(1\ts\theta^0,\dots,1\ts\theta^{m-1})$ of $L\ts_K L$, is in $T\ts_S T$, i.e.\ if and only if $\eta V_\tau^{-1}$ has entries in $T$. Which, in turn, is equivalent to 
$\eta (I - L_\tau) Y_\tau^{-1}$ having entries in $T$ (\ref{PropDiag7}).

Ad (iii). We use $E_\tau V_\tau = Y_\tau M_\tau$ (\ref{LemDiag6}).
\end{Proposition}

\subsection{A projective resolution of $T$ over $T\ts_S T$}
\label{SubSecProjRes}

\begin{quote}
\begin{footnotesize}
In the introduction to [LL 92], several sources for a projective resolution of $T$ over $T\ts_S T$ are indicated. We give still another alternative way to view 
such a projective resolution, using our isomorphic copy $(T\ts_S T)\delr_{T/S}$. We include the case of Hochschild (co)homology with galois twisted coefficients, for there a shift in the $2$-periodic 
vanishing of these groups occurs when compared to the untwisted case.
\end{footnotesize}
\end{quote}

\begin{Notation}
\label{NotPR1}\rm
We denote $\Lambda := (T\ts_S T)\delr_{T/S}$, i.e.\
\[
T\ts_S T\;\lraisoa{\delr_{T/S}}\;\Lambda\;\hra\;\prod_{\sigma\in G} T.
\]
For $i\in [0,m-1]$, we dispose of the projection map 
\[
\begin{array}{rcl}
\Lambda                 & \lraa{\pi_i}     & T \\
(\eta_j)_{j\in [0,m-1]} & \lramaps         & \eta_i, \\
\end{array}
\]
which is a ring morphism, and by means of which $T$ becomes a module over $\Lambda$, denoted by $T_i$. 

If we identify along $\delr_{T/S}$, the module operation of $x\ts y\in T\ts_S T$ on $z\in T_i$ is given by $z \cdot (x\ts y) := xzy^{\sigma_i}$. I.e.\ $T_i$ may be viewed as $T$ 
equipped with a structure as a galois twisted $T$-$T$-bimodule. In particular, $T_0$ may be viewed as $T$ equipped with the structure as an untwisted $T$-$T$-bimodule, that is 
$z \cdot (x\ts y) := xzy$. Note that $\pi_0$ is just the multiplication map.

In $\Lambda$, we have the elements 
\[
\begin{array}{rcl}
a & := & (\prod_{k\in [1,m-1]}(\theta - \theta^{\sigma_k}),0,\dots,0) \\
b & := & (0,\theta^{\sigma_1} - \theta,\dots, \theta^{\sigma_{m-1}} - \theta). \\
\end{array}
\]
The element $a$ is in $\Lambda$ by (\ref{ThEld4} ii), and $b$ is in $\Lambda$ by \mb{(\ref{ThEld4} iii)}. The multiplication map by $a$ resp.\ by $b$ shall be denoted by $\Lambda\lraa{\alpha}\Lambda$ resp.\ by 
$\Lambda\lraa{\beta}\Lambda$.
\end{Notation}

\begin{Lemma}
\label{PropPR2}
We have a $2$-periodic projective resolution
\[
\cdots\lraa{\alpha}\Lambda \lraa{\beta}\Lambda \lraa{\alpha}\Lambda \lraa{\beta}\Lambda \lraa{\pi_0} T_0 
\]
of $T_0$ over $\Lambda$.

\rm
From $ab = 0$ we take $\alpha\beta = 0$ and $\beta\alpha = 0$. Moreover, $\beta\pi_0 = 0$. By (\ref{ThEld4} ii), we obtain that the image of $\alpha$ equals the kernel of $\beta$. 
To prove that the image of $\beta$ equals the kernel of $\alpha$ (resp.\ of $\pi_0$), means, after identification along $\delr_{T/S}$, to show that this kernel is $T\ts_S T$-linearly 
generated by $b = 1\ts\theta - \theta\ts 1$. In fact, suppose given $\sum_{i\in [0,m-1]} u_i\ts\theta^i$ such that $\sum_{i\in [0,m-1]} u_i\theta^i = 0$, we may write 
$\sum_{i\in [0,m-1]} u_i\ts\theta^i = \sum_{i\in [0,m-1]} u_i (1\ts\theta^i - \theta^i\ts 1)$, and $1\ts\theta^i - \theta^i\ts 1$ is a multiple of $1\ts\theta - \theta\ts 1$.
\end{Lemma}

\begin{Proposition}[{cf.\ [BACH 91, prop.\ 1.3], [LL 92, 1.6.2]}]
\label{CorPR3}
\Absatz
Let $\phi := \sum_{k\in [1,m-1]} v_t(\theta - \theta^{\sigma_k})$. The Hochschild homology of $T$, over the ground ring $S$ and with coefficients in $T_i$, $i\in [0,m-1]$, is given in dimension 
$j\geq 0$ by
\[
H_j(T,T_0;S) \;\iso\;
\left\{
\begin{array}{ll}
T           & \mb{for $j = 0$} \\
T/t^{\phi}T & \mb{for $j$ odd} \\
0           & \mb{for $j$ even, $j\geq 2$} \\
\end{array}
\right.
\]
and by 
\[
H_j(T,T_i;S) \;\iso\;
\left\{
\begin{array}{ll}
0                                & \mb{for $j$ odd} \\
T/(\theta^{\sigma_i} - \theta) T & \mb{for $j$ even} \\
\end{array}
\right.
\]
for $i\in [1,m-1]$. The Hochschild cohomology is given by
\[
H^j(T,T_0;S) \;\iso\;
\left\{
\begin{array}{ll}
T           & \mb{for $j = 0$} \\
0           & \mb{for $j$ odd} \\
T/t^{\phi}T & \mb{for $j$ even, $j\geq 2$} \\
\end{array}
\right.
\]
and by 
\[
H^j(T,T_i;S) \;\iso\;
\left\{
\begin{array}{ll}
T/(\theta^{\sigma_i} - \theta) T & \mb{for $j$ odd} \\
0                                & \mb{for $j$ even} \\
\end{array}
\right.
\]
for $i\in [1,m-1]$.

\rm
Tensoring the resolution in (\ref{PropPR2}) with $T_i$ over $\Lambda$ yields the homology. Application of $\Hom_\Lambda(-,T_i)$ yields the cohomology.
\end{Proposition}

\subsection{The local ring $T\ts_S T$}

Assume $\theta = t$. 

\begin{Remark}
\label{RemLR1}\rm
The radical of $\Lambda$ is given by $\rfk\Lambda = \Lambda\cap \prod_{i\in [0,m-1]} tT$. In fact, the $(\phi+1)$st power of this intersection is contained in $t\Lambda$ (\ref{ThEld4} ii). And 
conversely, by (\ref{ThEld4} iii), this intersection has the $T$-linear basis
\[
\Big\{ (t,\dots,t)\Big\} \cup \left\{\left({\text\prod_{k\in [0,i-1]}} (t^{\sigma_j} - t^{\sigma_k})\right)_{j\in [0,m-1]} \;\Big|\; i\in [1,m-1] \right\},
\]
whence  $\Lambda/(\Lambda\cap \prod_{i\in [0,m-1]} tT)$ is isomorphic to $T_0/tT_0$. In particular, $T\ts_S T$ is a local ring.

Identifying along $\delta_{T/S}$, a $T$-linear basis of the radical $\rfk\Lambda$ is given by 
\[
(t\ts 1,1\ts t,1\ts t^2,\dots,1\ts t^{m-1}). 
\]
Given $j\in\Z$, we write $\ul{j} := \max(j,0)$. Suppose given $i\geq 0$. A basis of $\rfk^i\Lambda$ is given by
\[
(t^{\ul{i-0}}\ts 1,t^{\ul{i-1}}\ts t,t^{\ul{i-2}}\ts t^2,\dots,t^{\ul{i-(m-1)}}\ts t^{m-1}).
\]
In particular, we have
\[
\dim_{T/tT} \rfk^i\Lambda/\rfk^{i+1}\Lambda = \min(i+1,m).
\]
Cf.\ [K\"u 99, E.2.3].
\end{Remark}

\section{The cyclotomic Dedekind embedding}
\label{SecCyc}

\begin{quote}
\begin{footnotesize}
We shall apply (\ref{ThEld4}) to the case of the extension $\Z_{(p)}\tm\Z_{(p)}[\zeta_{p^n}]$.
\end{footnotesize}
\end{quote}

\begin{Setup}
\label{NotVand1}\rm
Let $p$ be a prime, let $n\geq 1$. In the notation of (\ref{SetupEld3}), we place ourselves in the situation $S = \Z_{(p)}$, $s = p$, $T = \Z_{(p)}[\zeta_{p^n}]$ 
and $t = \theta = 1 - \zeta_{p^n}$.

We consider the Dedekind embedding
\[
\begin{array}{rclcl}
\Z_{(p)}[\zeta_{p^n}] & \ts & \Z_{(p)}[\zeta_{p^n}] & \hrah{\delta_{p^n}} & \prod_{j\in(\sZ/p^n)^\ast}\Z_{(p)}[\zeta_{p^n}] \\
\zeta_{p^n}^k         & \ts & \zeta_{p^n}^l         & \lramaps            & (\zeta_{p^n}^k\zeta_{p^n}^{jl})_{j\in(\sZ/p^n)^\ast} \\
\end{array}
\]
where $k,l\in [0,p^n-1]$. With respect to the $\Z_{(p)}[\zeta_{p^n}]$-linear basis $(1\ts t^i)_{i\in [0,(p-1)p^{n-1}-1]}$ of $\Z_{(p)}[\zeta_{p^n}] \ts \Z_{(p)}[\zeta_{p^n}]$ and to the 
tuple basis on the right hand side, this embedding is given by the Vandermonde matrix $V_\tau$ for $\tau := (1-\zeta_{p^n}^j)_{j\in [0,p^n-1]\ohne (p)}$.
\end{Setup}

\begin{Lemma}
\label{RemVand2}
For $i,j\in [0,p^n-1]$ with $i\neq j$, we have $v_t(\zeta_{p^n}^i - \zeta_{p^n}^j) = (i-j)[p]$.

\rm
We may assume $j = 0$ and $i = i[p]$. Using $v_t(p) = (p-1)p^{n-1}$ (resp.\ using a direct calculation if $p = 2$, $i = 2^{n-1}$), the congruence $\zeta_{p^n}^i - 1\con_p (\zeta_{p^n} - 1)^i$ 
yields $v_t(\zeta_{p^n}^i - 1) = v_t((\zeta_{p^n} - 1)^i) = i$.
\end{Lemma}

\begin{Lemma}
\label{RemVand3}
Suppose given $j\geq 0$, and write it $p$-adically as $j = \sum_{l\geq 0} a_l p^l$, where $a_l\in [0,p-1]$. Then
\[
\sumd{i\in [1,j]} i[p] \; =\; \sumd{k\geq 0} (a_k - a_{k+1})(k+1)p^k.
\]

\rm
Denote left and right hand side of the claimed equation by $l(j)$ and $r(j)$, respectively. We have $l(0) = 0 = r(0)$. Moreover, for $j\geq 1$ we have $l(j) - l(j-1) = j[p]$.
If $v_p(j) = 0$, then $r(j) - r(j-1) = 1$, so that we may suppose $v_p(j)\geq 1$. Writing $j-1 = \sum_{k\geq 0} a'_k p^k$, $a'_k\in [0,p-1]$, we note that $a'_k = p-1$ and $a_k = 0$ 
for $k\leq v_p(j) - 1$, $a'_k = a_k-1$ for $k = v_p(j)$ and $a'_k = a_k$ for $k \geq v_p(j)+1$. Therefore, 
\[
\begin{array}{rcl}
r(j) - r(j-1)
& = & \sumd{k\geq 0} (a_k - a'_k)(k+1)p^k - \sumd{k\geq 1} (a_k - a'_k) k p^{k-1}\\
& = & (v_p(j) + 1) j[p]\; - \sumd{k\in [0,v_p(j) - 1]} (p-1)(k+1)p^k \\
& - & v_p(j) j[p]/p\; + \sumd{k\in [0,v_p(j) - 2]} (p-1)(k+1)p^k \\
& = & j[p]. \\
\end{array}
\]
\end{Lemma}

\begin{Lemma}
\label{RemVand7}
Suppose given $k\geq j\geq 0$, with $j, k\not\in (p)$. Write $j = \sum_{l\geq 0} a_l p^l$, $k = \sum_{l\geq 0} a'_l p^l$, $k-j = \sum_{l\geq 0} a''_l p^l$, where $a_l,a'_l,a''_l\in [0,p-1]$.
Then
\[
\sumd{i\in [0,j-1]\ohne (p)} (k-i)[p] = -1 + \sumd{l\geq 0} \left(\left((a'_l - a''_l) - (a'_{l+1} - a''_{l+1})\right)(l+1) - a_{l+1}\right) p^l.
\]
In particular,
\[
\sumd{i\in [0,j-1]\ohne (p)} (j-i)[p] = -1 + \sumd{l\geq 0} \left(a_l(l+1) - a_{l+1}(l+2)\right) p^l.
\]

\rm
This follows by (\ref{RemVand3}) and the remark that $\sumd{i\in [0,j-1]\cap (p)} (k-i)[p] = \#([0,j-1]\cap (p)) = 1 + \sumd{l\geq 0} a_{l+1}p^l$.
\end{Lemma}

\begin{Lemma}
\label{LemVand7_5}
Keep the notation of (\ref{RemVand7}), but allow $k\geq j\geq 0$ to be arbitrary. We have
\[
\sumd{l\geq 0} \left((a'_l - a''_l - a_l) - (a'_{l+1} - a''_{l+1} - a_{l+1})\right)(l+1) p^l \geq 0.
\]

\rm
Let $U := \{ l\geq 0\; |\; a_l + a''_l\geq p\} \tm \Z_{\geq 0}$. We obtain
\[
\begin{array}{rl}
  & \sumd{l\geq 0} \left((a'_l - a''_l - a_l) - (a'_{l+1} - a''_{l+1} - a_{l+1})\right)(l+1)p^l \\
= & \sumd{l\in U} (-p)(l+1)p^l + \sumd{l\in U+1} (l+1)p^l -  \sumd{l\in U-1} (-p)(l+1)p^l - \sumd{l\in U} (l+1)p^l \\
= & \sumd{l\in U} \left(-(l+1)p^{l+1} + (l+2) p^{l+1} + lp^l - (l+1)p^l\right) \\
= & \sumd{l\in U} (p - 1)p^l. \\
\end{array}
\]
\end{Lemma}

\begin{Lemma}
\label{LemVand8}
The tuple $\tau = (1-\zeta_{p^n}^j)_{j\in [0,p^n-1]\ohne (p)}$ is minimally ordered. 

\rm
Using (\ref{RemVand2}), we need to see that 
\[
\sumd{i\in [0,j-1]\ohne (p)} (j-i)[p] \leq \sumd{i\in [0,j-1]\ohne (p)} (k-i)[p]
\]
for $j\in [0,p^n-1]\ohne (p)$ and $k\in [j+1,p^n-1]\ohne (p)$. In the notation and using the assertion of (\ref{RemVand7}) this amounts to the inequality
\[
\sumd{l\geq 0} \left((a'_l - a''_l - a_l) - (a'_{l+1} - a''_{l+1} - a_{l+1})\right)(l+1) p^l \geq 0
\]
treated in (\ref{LemVand7_5}).
\end{Lemma}

\begin{Theorem}
\label{PropVand9}
Suppose given $j\in [0,p^n - 1]\ohne (p)$. Let $N(j) := \#\left([0,j]\ohne (p)\right)$. Write $j = \sum_{s\geq 0} a_s p^s$ with $a_s\in [0,p-1]$. The $N(j)$th 
$\Z_{(p)}[\zeta_{p^n}]$-linear elementary divisor of the cyclotomic Dedekind embedding 
\[
\begin{array}{rclcl}
\Z_{(p)}[\zeta_{p^n}] & \ts & \Z_{(p)}[\zeta_{p^n}] & \hrah{\delta_{p^n}} & \prod_{i\in(\sZ/p^n)^\ast}\Z_{(p)}[\zeta_{p^n}] \\
\zeta_{p^n}^k         & \ts & \zeta_{p^n}^l         & \lramaps            & (\zeta_{p^n}^{k+il})_{i\in(\sZ/p^n)^\ast}, \\
\end{array}
\]
where $k,l\in [0,(p-1)p^{n-1}-1]$, has valuation
\[
-1 + \sumd{s\geq 0} \left(a_s(s+1) - a_{s+1}(s+2)\right) p^s
\]
at $t$.

\rm
This follows by (\ref{ThEld4} i) using (\ref{LemVand8}, \ref{RemVand7}, \ref{RemVand2}).
\end{Theorem}

\begin{quote}
\begin{footnotesize}
\begin{Remark}
\label{RemVand9_5}\rm
{\sc Plesken} gives a system of ties that describes the image of the Dedekind embedding case $n = 1$ [P 80, p.\ 60].
\end{Remark}

\begin{Remark}
\label{RemVand10}
\rm
If $n = 1$, the valuation at $t$ of the determinant of this embedding $\delta_p$ is $-(p-1) + \sum_{a\in [1,p-1]} a = (p-1)(p-2)/2$. If $n\geq 2$, we obtain the valuation at $t$ of the 
determinant of $\delta_{p^n}$ to be
\[
\begin{array}{rlcl}
  & -(p-1)p^{n-1} & + & \sumd{(a_l)_{l\in [0,n-1]}\in [0,p-1]^n,\; a_n \; =\; 0,\; a_0\neq 0}\;\;\;\;\sumd{l\in [0,n-1]} (a_l(l+1) - a_{l+1}(l+2))p^l \\
= & -(p-1)p^{n-1} & + & p^{n-2}\cdot \sumd{a\in [1,p-1],\; b\in [0,p-1]} (a - 2b) \\
  &               & + & \sumd{l\in [1,n-2]} (p-1)p^{n-3}\cdot \sumd{a,b\in [0,p-1]} (a(l+1) - b(l+2))p^l \\
  &               & + & (p-1)p^{n-2}\cdot \sumd{a\in [0,p-1]} a n p^{n-1}  \\
= & -(p-1)p^{n-1} & + & p^{n-2}\cdot (p-1)(-p^2 + 2 p)/2 \\
  &               & - & (p-1)p^{n-2}\cdot p(p-1)/2 \cdot p(p^{n-2}-1)/(p-1) \\
  &               & + & (p-1)p^{n-2}\cdot p(p-1)/2 \cdot n p^{n-1} \\
= &               &   & \hspace*{-27mm} p^{2n-2}(p-1)((p-1)n - 1)/2. \\
\end{array}
\]
Hence for $n\geq 1$, we recalculated the valuation at $p$ of the discriminant of $\Z_{(p)}[\zeta_{p^n}]$ over $\Z_{(p)}$ to be $p^{n-1}((p-1)n - 1)$. Cf.\ [N 91, 10.1] (or \ref{RemCP9_0_6}).
\end{Remark}

\begin{Example}
\label{ExVand11}\rm
Let $p = 3$ and $n = 2$, so that
\[
\tau \; =\; (1-\zeta_9,\, 1-\zeta_9^2,\, 1-\zeta_9^4,\, 1-\zeta_9^5,\, 1-\zeta_9^7,\, 1-\zeta_9^8).
\]
As elementary divisors of the Dedekind embedding
\[
\begin{array}{rclcl}
\Z_{(3)}[\zeta_9] & \ts & \Z_{(3)}[\zeta_9] & \hraa{\delta_9} & \prod_{j\in(\sZ/9)^\ast}\Z_{(3)}[\zeta_{9}] \\
\zeta_9^k         & \ts & \zeta_9^l         & \lramaps        & (\zeta_9^k\zeta_9^{jl})_{j\in(\sZ/9)^\ast}, \\
\end{array}
\]
where $k,l\in [0,5]$, (\ref{PropVand9}) yields
\[
(t^0,\,t^1,\,t^4,\,t^5,\,t^8,\,t^9).
\]
The determinant of $\delr_9$ has valuation $27$ at $t$ (cf.\ \ref{RemVand10}). By (\ref{ThEld4} iii), a triangular $\Z_{(3)}[\zeta_9]$-linear basis of the image of $\delr_9$ is given by
the rows of the matrix
\[
\left[
\begin{array}{cccccc}
\scm 1 &\scm 1 &\scm 1 &\scm 1 &\scm 1 &\scm 1 \\
\scm 0 &\scm - \zeta_9 + \zeta_9^2 &\scm - \zeta_9 + \zeta_9^3 &\scm - \zeta_9 + \zeta_9^4 &\scm -2\zeta_9 - \zeta_9^4 &\scm - \zeta_9 - \zeta_9^2 - \zeta_9^5 \\
\scm 0 &\scm 0 &\scm -1 - \zeta_9^4 - \zeta_9^5  &\scm 1 - \zeta_9^2 + 2\zeta_9^3 - 2\zeta_9^5 &\scm -1 + \zeta_9^2 + \zeta_9^3 + 2\zeta_9^5 &\scm -1 - 2\zeta_9 + \zeta_9^3 - \zeta_9^4 \\
\scm 0 &\scm 0 &\scm 0 &\scm 1 - 2\zeta_9 - 2\zeta_9^2 + 2\zeta_9^3 -\zeta_9^4 - \zeta_9^5 &\scm 2 + 2\zeta_9 + 2\zeta_9^2 + 4\zeta_9^3 + \zeta_9^5 
  &\scm -1 - \zeta_9 + 2\zeta_9^2 + \zeta_9^3 + \zeta_9^4 + \zeta_9^5 \\
\scm 0 &\scm 0 &\scm 0 &\scm 0 &\scm 3 + 6\zeta_9 + 3\zeta_9^3 - 3\zeta_9^5 &\scm 3\zeta_9 + 3\zeta_9^2 + 3\zeta_9^3 + 3\zeta_9^4 + 3\zeta_9^5 \\
\scm 0 &\scm 0 &\scm 0 &\scm 0 &\scm 0 &\scm 3\zeta_9^2 + 3\zeta_9^4 + 3\zeta_9^5 \\
\end{array}
\right].
\]
\end{Example}
\end{footnotesize}
\end{quote}
\section{The cyclic Wedderburn embedding}

\subsection{The q-Pascal method}

\begin{quote}
\begin{footnotesize}
We diagonalize the Vandermonde matrix that describes the Wedderburn embedding in the $\Z[\zeta_m]$-linear basis $(c_m^0,\dots,c_m^{m-1})$ of $\Z[\zeta_m] C_m$. More specifically speaking, we multiply 
with a lower triangular matrix $G_q$ that records the $q$-Pascal triangle, consisting of Gau\ss ian polynomials, from the right and with its transpose $G_q^\trp$ from the left. The resulting diagonal 
matrix is of elementary divisor form, regardless whether or not the ground ring $\Z[\zeta_m]$ is a principal ideal domain. We derive the necessary identities firstly in a pedestrian fashion, 
yielding in particular a formula for $G_q^{-1}$, and secondly, as a consequence of {\sc F.\ H.\ Jackson}'s $q$-analogue (1910) of {\sc L.\ Saalsch\"utz'} theorem (1890).
\end{footnotesize}
\end{quote}

\begin{Notation}
\label{NotSaal0}\rm
Let $m\geq 1$. The {\it cyclic Wedderburn embedding} is given by
\[
\begin{array}{rcl}
\Z[\zeta_m]C_m & \hraa{\omega_m} & \prod_{j\in [0,m-1]} \Z[\zeta_m] \\
c_m            & \lramaps        & (\zeta_m^j)_{j\in [0,m-1]}. \\
\end{array}
\]
\end{Notation}

\subsubsection{Pedestrian}

\begin{Notation}[Gau\ss ian polynomials]
\label{NotSaal1}\rm
Consider the field $\Q(q)$ of rational functions in the indeterminate $q$. Given $i\geq 0$ and $j\in\Z$ and another indeterminate $t$, we let
\[
\begin{array}{rcl}
{[i]}            & := & (q^i-1)/(q-1) \\
{[i]}!           & := & \prod_{k\in [1,i]} [k]\vspace*{1mm} \\
\smateckze{i}{j} & := & 
\left\{
\begin{array}{ll}
\frac{[i]!}{[j]![i-j]!} & \mb{if $j\in [0,i]$} \\
0                        & \mb{if $j\not\in [0,i]$} \\
\end{array}
\right. \\
(t;q)_i          & := & \prod_{k\in [0,i-1]} (1 - q^k t). \\
\end{array}
\]
In particular, $[0]! = 1$. If necessary, we indicate $q$ as $\smateckze{i}{j}_q := \smateckze{i}{j}$. Let the {\it $q$-Pascal matrix} of size $m\ti m$ be defined by
\[
G_q \; :=\; \left( \smateckze{i}{j} \right)_{i,j\in [0,m-1]}.
\]
Moreover, we let
\[
\begin{array}{rcl}
I   & := & \left(\dell_{i,j}\right)_{i,j\in [0,m-1]} \\
V_q & := & \left(q^{ij}\right)_{i,j\in [0,m-1]} \\
D_q & := & \left(\dell_{i,j} [i]! (q-1)^i q^{\ssmatze{i}{2}}\right)_{i,j\in [0,m-1]}. \\
\end{array}
\]
\end{Notation}

\begin{Lemma}
\label{LemSaal2}
If $i\geq 1$ and $j\geq 0$, we have
\[
\begin{array}{rcrcrl}
\smateckze{i}{j} & = &         \smateckze{i-1}{j-1} & + & q^j \smateckze{i-1}{j} & \\
                 & = & q^{i-j} \smateckze{i-1}{j-1} & + &     \smateckze{i-1}{j} & . \\
\end{array}
\]
Moreover,
\[
(t;q)_i \; = \; \sumd{k\in [0,i]} (-1)^k q^{\ssmatze{k}{2}} \smateckze{i}{k} t^k. 
\]
Finally,
\[
\smateckze{i}{j}_{q^{-1}} \; = \; q^{-j(i-j)}\smateckze{i}{j}_q.
\]
\end{Lemma}

\begin{Proposition}[Inversion of the $q$-Pascal matrix]
\label{PropSaal3}
\Absatz
We obtain the inverse of $G_q$ to be
\[
G_q^{-1} = \left( (-1)^{j+k}q^{\ssmatze{j-k}{2}}\smateckze{j}{k} \right)_{j,k\in [0,m-1]}.
\]

\rm
Given $0\leq k\leq i\leq m-1$, we need to show that
\[
\sumd{j\in [k,i]} \smateckze{i}{j}\smateckze{j}{k} (-1)^{j+k} q^{\ssmatze{j-k}{2}} \; = \; \dell_{i,k}.
\]
We perform an induction on $i-k$. If $i = k$, we obtain $1 = 1$. If $i = k+1$, we obtain $0 = 0$. If $i \geq k + 2$, we obtain
\[
\begin{array}{rcl}
\sumd{j\geq 0} \smateckze{i}{j}\smateckze{j}{k} (-1)^{j+k} q^{\ssmatze{j-k}{2}} 
& \aufgl{(\ref{LemSaal2})} & \sumd{j\geq 0} \smateckze{i-1}{j}\smateckze{j}{k} (-1)^{j+k} q^{\ssmatze{j-k}{2}} \\
& + & \sumd{j\geq 0} \smateckze{i-1}{j-1}\smateckze{j-1}{k} (-1)^{j+k} q^{\ssmatze{j-k}{2} + (i-j)} \\
& + & \sumd{j\geq 0} \smateckze{i-1}{j-1}\smateckze{j-1}{k-1} (-1)^{j+k} q^{\ssmatze{j-k}{2} + (i-j) + (j-k)} \\
& \aufgl{induction} & \sumd{j\geq 0}\smateckze{i-1}{j-1}\smateckze{j-1}{k-1} (-1)^{j+k} q^{\ssmatze{j-k}{2} + i - k}. \\ 
\end{array}
\]
The last term vanishes by subinduction on $k$, since if $k = 0$ we obtain
\[
\sumd{j\geq 0} \smateckze{i}{j} (-1)^j q^{\ssmatze{j}{2}} 
\;\aufgl{(\ref{LemSaal2})}\; (1;q)_i \; =\;  0.
\]
\end{Proposition}

\begin{Lemma}
\label{PropSaal4}
We have
\[
V_q (G_q^{-1})^\trp \; = \; G_q D_q.
\]

\rm
We need to show that for $i,k\in [0,m-1]$ 
\[
\sumd{j\in [0,k]} q^{ij} (-1)^{j+k} \smateckze{k}{j} q^{\ssmatze{k-j}{2}} \; =\; q^{\ssmatze{k}{2}} [k]!(q-1)^k \smateckze{i}{k}.
\]
We perform an induction on $k$, starting in case $k = 0$ with $1 = 1$. If $k\geq 1$, we calculate
\[
\begin{array}{rcl}
\sumd{j\in [0,k]} q^{ij} (-1)^{j+k} \smateckze{k}{j} q^{\ssmatze{k-j}{2}}
& \aufgl{(\ref{LemSaal2})} & \sumd{j\in [0,k-1]} q^{i(j+1)} (-1)^{j+1+k} \smateckze{k-1}{j} q^{\ssmatze{k-j-1}{2}} \\
& + & \sumd{j\in [0,k-1]} q^{ij} (-1)^{j+k} \smateckze{k-1}{j} q^{\ssmatze{k-j}{2}+j} \\
& \aufgl{induction} & (q^i - q^{k-1}) q^{\ssmatze{k-1}{2}} [k-1]!(q-1)^{k-1} \smateckze{i}{k-1} \\
& = & q^{\ssmatze{k}{2}} [k]!(q-1)^k \smateckze{i}{k}.
\end{array}
\]
\end{Lemma}

\begin{Lemma}[Fourier inversion]
\label{LemIF1}
We have $V_{\zeta_m}V_{\zeta_m^{-1}} = mI$. In particular, reordering the rows of $V_{\zeta_m^{-1}}$ shows that
\[
(\det V_{\zeta_m})^2 \; =\;
\left\{
\begin{array}{ll}
(-1)^{\frac{m-1}{2}} m^m & \mb{if $m$ is odd} \\
(-1)^{\frac{m-2}{2}} m^m & \mb{if $m$ is even.} \\
\end{array}
\right.
\]

\rm
In fact,
\[
\begin{array}{rcl}
\sum_{j\in [0,m-1]} \zeta_m^{ij}\zeta_m^{-jk}
& = & \sum_{j\in [0,m-1]} \zeta_m^{(i-k)j} \\
& = & m\dell_{i,k}. \\
\end{array}
\]
\end{Lemma}

\begin{quote}
\begin{footnotesize}
\begin{Corollary}[{cf.\ [N 91, 8.6]}]
\label{CorIF1_5}
If $m$ is odd, then $\Q(\zeta_m)$ contains $\sqrt{(-1)^{\frac{m-1}{2}}m}$. 
\end{Corollary}
\end{footnotesize}
\end{quote}

\begin{Proposition}
\label{ThSaal5}
\Absit
\begin{itemize}
\item[(i)] We have the diagonalization
\[
G_{\zeta_m^{-1}}^\trp V_{\zeta_m} G_{\zeta_m^{-1}} = m (D_{\zeta_m^{-1}})^{-1},
\]
the right hand side being a diagonal matrix with $(i+1)$st diagonal entry
\[
\fracd{m \zeta_m^{(i^2)}}{\prod_{j\in [1,i]} (1-\zeta_m^j)},
\]
where $i\in [0,m-1]$.

\item[(ii)] We have
\[
\begin{array}{rl}
    & (\Z[\zeta_m] C_m)\omega_m \\
=   & \left\{ (y_j)_{j\in [0,m-1]}\;\left|\rule[7mm]{0mm}{0mm}\right.\; 
\sumd{j\in [i,m-1]} y_j \smateckze{j}{i}_{\zeta_m^{-1}} \in \Z[\zeta_m]\cdot\fracd{m}{\prod_{j\in [1,i]} (1-\zeta_m^j)} \;\mb{\rm\ for $i\in [0,m-1]$} \right\} \\
\tm & \prod_{j\in [0,m-1]}\Z[\zeta_m].
\end{array}
\]

\item[(iii)] A $\Z[\zeta_m]$-linear basis of $(\Z[\zeta_m] C_m)\omega_m$ is given by the following tuple of elements of $\prod_{j\in [0,m-1]}\Z[\zeta_m]$.
\[
\left(\left((-1)^k \zeta_m^{\ssmatze{k}{2}}{\text\fracd{m}{\prod_{l\in [1,j]}(1 - \zeta_m^l)}} \smateckze{j}{k}_{\zeta_m}\right)_{\! k\in [0,m-1]} \right)_{\! j\in [0,m-1]}
\]
\end{itemize}

\rm
Ad (i). This formula follows from (\ref{PropSaal4}, \ref{LemIF1}).

Ad (ii). A vector $y = (y_j)_{j\in [0,m-1]}$ is contained in the image of $\omega_m$ if and only if $y V_{\zeta_m}^{-1}$ has entries in $\Z[\zeta_m]$. But 
$V_{\zeta_m}^{-1} (G_{\zeta_m^{-1}}^\trp)^{-1} = G_{\zeta_m^{-1}} D_{\zeta_m^{-1}} m^{-1}$, so this condition translates into
$y G_{\zeta_m^{-1}} D_{\zeta_m^{-1}} m^{-1}$ to have entries in $\Z[\zeta_m]$, which is the defining condition given above. 

Ad (iii). This follows by $G_{\zeta_m^{-1}}^\trp V_{\zeta_m} = m (D_{\zeta_m^{-1}})^{-1} (G_{\zeta_m^{-1}})^{-1}$ using (\ref{PropSaal3}, \ref{LemSaal2}).
\end{Proposition}

\begin{Remark}
\label{RemSaal7}\rm
Because the derivative of
\[
F(X) \erstgl X^m - 1 \zweigl \prod_{j\in [0,m-1]} (X - \zeta_m^j)
\]
evaluated at $1$ yields
\[
F'(1) \erstgl m \zweigl\prod_{j\in [1,m-1]} (1 - \zeta_m^j),
\]
the $m$th diagonal entry in (\ref{ThSaal5} i) equals $\zeta_m$. This reproves in particular that all diagonal entries of (\ref{ThSaal5} i) are contained in $\Z[\zeta_m]$.
\end{Remark}

\begin{Remark}
\label{RemSaal8}\rm
Using (\ref{RemSaal7}), for $n\geq 1$, $p$ prime, $m = p^n$, the valuation at $t = 1 - \zeta_{p^n}$ of the $(i+1)$th elementary divisor of $\omega_{p^n}$, $i\in [0,p^n-1]$, is given by
\[
\begin{array}{rcl}
v_t\left(\fracd{p^n}{\prod_{j\in [1,p^n-1-i]} (1-\zeta_{p^n}^j)}\right) 
& = & \sum_{j\in [1,i]} v_t\left(1-\zeta_{p^n}^{-j}\right) \\
& \aufgl{(\ref{RemVand2})} & \sum_{j\in [1,i]} j[p], \\
\end{array}
\] 
in accordance with (\ref{PropVand5}) below in view of (\ref{RemVand3}).
\end{Remark}

\subsubsection{Invoking $q$-Saalsch\"utz}
\label{SubSecSaal}

Following {\sc J.\ Stokman,} we cite the

\begin{Theorem}[{{\sc Saalsch\"utz,} {\sc Jackson} [A 76, 3.3.12] {\rm (\footnotemark)}}]
\label{ThSaal9}
Given $n\geq 0$, and further indeterminates $a$, $b$, $c$, we have
\[
\sum_{m\geq 0} \fracd{(a;q)_m (b;q)_m (q^{-n};q)_m}{(q;q)_m (c;q)_m (abq^{1-n}c^{-1};q)_m}\cdot q^m \; = \; \frac{(c/a;q)_n (c/b;q)_n}{(c;q)_n (c/ab;q)_n}.
\]
\end{Theorem}
\footnotetext{\scr In the left hand side summand of [A 76, 3.3.14], a factor $(ab/c)^{N-n}$ has been forgotten.}

and specialize to the

\begin{Corollary}[{$\equ$ \ref{PropSaal4}}]
\label{CorSaal10}
Given $m\geq 1$, we have
\[
V_q  \; = \; G_q D_q (G_q)^\trp.
\]

\rm
Given $i,k\in [0,m-1]$, we need to show that
\[
q^{ik} = \sum_{j\in [0,m-1]} \smateckze{i}{j} \cdot [j]! (q-1)^j q^{\ssmatze{j}{2}} \cdot \smateckze{k}{j}. 
\]
Substituting $q$ by $q^{-1}$, $n$ by $k$ and specializing $c$ and then $b$ to zero (after rewriting the right hand side), and specializing $a$ to $q^i$, (\ref{ThSaal9}) becomes
\[
\begin{array}{rcl}
q^{ik} 
& = & \sum_{j\geq 0} \fracd{(q^i;q^{-1})_j (q^k;q^{-1})_j}{(q^{-1};q^{-1})_j}\cdot q^{-j} \vspace*{1mm}\\
& = & \sum_{j\geq 0} \fracd{(1-q^i)\cdots (1-q^{i-j+1})\cdot (1-q^k)\cdots (1-q^{k-j+1})}{(1-q^{-1})\cdots (1-q^{-j})}\cdot q^{-j} \vspace*{1mm}\\
& = & \sum_{j\geq 0} \fracd{(q^i-1)\cdots (q^{i-j+1}-1)}{(q^1-1)\cdots (q^j-1)} \cdot \fracd{(q^k-1)\cdots (q^{k-j+1}-1)}{(q^1-1)\cdots (q^j-1)} \cdot q^{\ssmatze{j}{2}} [j]! (q-1)^j \vspace*{1mm}\\
& = & \sum_{j\geq 0} \smateckze{i}{j} \cdot \smateckze{k}{j} \cdot q^{\ssmatze{j}{2}} [j]! (q-1)^j. \\
\end{array}
\]
\end{Corollary}

\subsection{The general Vandermonde method}

\begin{quote}
\begin{footnotesize}
To compare, and to recalculate elementary divisors, we apply the method of section \ref{SubSecElVand} to the localized cyclic Wedderburn embedding.
\end{footnotesize}
\end{quote}

\begin{Setup}
\label{SetupVand0}\rm
Let $p$ be a prime, let $n\geq 1$. In the notation of (\ref{SetupEld3}), we place ourselves in the situation $S = \Z_{(p)}$, $s = p$, $T = \Z_{(p)}[\zeta_{p^n}]$ 
and $t = \theta = 1 - \zeta_{p^n}$. We consider the localized cyclic Wedderburn embedding
\[
\begin{array}{rcl}
\Z_{(p)}[\zeta_{p^n}] C_{p^n} & \lraisoa{\omega_{p^n}} & \prod_{j\in\sZ/p^n} \Z_{(p)}[\zeta_{p^n}] \\
c_{p^n}^i                     & \lra                   & (\zeta_{p^n}^{ij})_{j\in\sZ/p^n}, \\
\end{array}
\]
where $i\in [0,p^n-1]$. In the notation of (\ref{NotDiag1}), and with respect to the $\Z_{(p)}[\zeta_{p^n}]$-linear basis $((1-c_{p^n})^i)_{i\in [0,p^n-1]}$ of 
$\Z_{(p)}[\zeta_{p^n}] C_{p^n}$ and the standard basis on the right hand side, this embedding is given by the Vandermonde matrix $V_\tau$ for $\tau := (1-\zeta_{p^n}^j)_{j\in [0,p^n-1]}$.
\end{Setup}

\begin{Lemma}
\label{LemVand4}
The tuple $\tau := (1-\zeta_{p^n}^j)_{j\in [0,p^n-1]}$ is minimally ordered.

\rm
Using (\ref{RemVand2}), we need to see that 
\[
\sumd{i\in [0,j-1]} (j-i)[p] \leq \sumd{i\in [0,j-1]} (k-i)[p]
\]
for $j\in [0,p^n-1]$ and $k\in [j+1,p^n-1]$. Write $j = \sum_{l\geq 0} a_l p^l$, $k = \sum_{l\geq 0} a'_l p^l$, $k-j = \sum_{l\geq 0} a''_l p^l$, where $a_l,a'_l,a''_l\in [0,p-1]$. By 
(\ref{RemVand3}), we need to show that
\[
\sumd{l\geq 0} \left((a'_l - a''_l - a_l) - (a'_{l+1} - a''_{l+1} - a_{l+1})\right)(l+1)p^l \geq 0. 
\]
This follows from (\ref{LemVand7_5}).
\end{Lemma}

\begin{Proposition}
\label{PropVand5}
Suppose given $j\in [0,p^n - 1]$. Write $j = \sum_{k\geq 0} a_k p^k$, $a_k\in [0,p-1]$. The $(j+1)$st elementary divisor of the localized cyclic Wedderburn embedding
\[
\begin{array}{rcl}
\Z_{(p)}[\zeta_{p^n}] C_{p^n} & \hrah{\omega_{p^n}} & \prod_{j\in\sZ/p^n} \Z_{(p)}[\zeta_{p^n}] \\
c_{p^n}^i                     & \lramaps            & (\zeta_{p^n}^{ij})_{j\in\sZ/p^n} \\
\end{array}
\]
where $i\in [0,p^n-1]$, has valuation
\[
\sumd{k\geq 0} (a_k - a_{k+1})(k+1)p^k
\]
at $t$.

\rm
This follows by (\ref{PropDiag7}, \ref{LemEld2}), using (\ref{LemVand4}, \ref{RemVand3}), or by (\ref{RemSaal8}), using (\ref{RemVand3}).
\end{Proposition}

\begin{quote}
\begin{footnotesize}
\begin{Remark}
\label{RemVand6}
\rm
By (\ref{LemIF1}), the square of the determinant of $\omega_{p^n}$ has valuation $np^n$ at $p$. Alternatively, by (\ref{PropVand5}), we obtain the valuation at $t$ of the determinant of this 
embedding to be
\[
\begin{array}{rl}
  & \sumd{(a_l)_{l\in [0,n-1]}\in [0,p-1]^n,\; a_n \; =\; 0}\;\;\;\;\sumd{l\in [0,n-1]} (a_l - a_{l+1})(l+1)p^l \\
= & \left(\sumd{l\in [0,n-2]} p^{n-2}\cdot \sumd{a,b\in [0,p-1]} (a - b)(l+1)p^l\right) + \left(p^{n-1}\cdot \sumd{a\in [0,p-1]} a n p^{n-1} \right)\\
= & p^{n-1} \frac{p(p-1)}{2} n p^{n-1}. \\
\end{array}
\]
\end{Remark}
\end{footnotesize}
\end{quote}
\subsection{The Pascal method}

\subsubsection{First order Pascal ties}
\label{SubSecPasWedFirst}

\begin{quote}
\begin{footnotesize}
There are some obvious ties, i.e.\ congruences of tuple entries, that are necessary for elements of $\prod_{i\in [0,p^n-1]}\Z[\zeta_{p^n}]$ to lie in $(\Z[\zeta_{p^n}]C_{p^n})\omega_{p^n}$.
If $n = 1$, they are already sufficient, yielding a manageable basis of $(\Z[\zeta_p]C_p)\omega_p$. In general, they describe an intermediate ring between $(\Z[\zeta_{p^n}]C_{p^n})\omega_{p^n}$ 
and $\prod_{i\in [0,p^n-1]}\Z[\zeta_{p^n}]$.
\end{footnotesize}
\end{quote}

\begin{Notation}
\label{NotPT0}
\rm
Suppose given $s\geq 0$, a polynomial $f(X) := \sum_{i\geq 0}a_i X^i\in\C[X]$ and a tuple $(y_j)_{j\in [0,s]}$, $y_j\in\C$. We define the {\it evaluation of $f$ at $(y_j)_{j\in [0,s]}$} to be
\[
\left(\sum_{i\geq 0}a_i X^i\right)\left[(y_j)_{j\in [0,s]}\right] := f\left[(y_j)_{j\in [0,s]}\right] := \sum_{j\in [0,s]} a_j y_j.
\]
We note the difference between the polynomial power $f^i\left[(y_j)_{j\in [0,s]}\right]$ (power of $f$, taken in $\C[X]$, evaluated) and the ordinary power $f\left[(y_j)_{j\in [0,s]}\right]^i$ 
(power of the evaluation of $f$, taken in $\C$), $i\geq 0$. 

Let $m\geq 1$, let $t = 1 - \zeta_m$. We consider the $\Z[\zeta_m]$-submodule $W^{(1)}_m$ of 
\[
\fbox{$\;\;
W^{(0)}_m := \prod_{j\in [0,m-1]}\Z[\zeta_m]
\;\;$}
\]
defined by the {\it first order Pascal ties}
\[
\fbox{$\;\;
\begin{array}{rcl}
\rule[7mm]{0mm}{0mm} W^{(1)}_m 
& := & \Big\{ (y_j)_{j\in [0,m-1]}\;\Big|\; \Big((1-X)^i\Big)\left[(y_j)_{j\in [0,m-1]}\right]\con_{t^i} 0 \mb{ for all } i\in [0,m-1] \Big\} \vspace*{2mm}\\
& =  & \Big\{ (y_j)_{j\in [0,m-1]}\;\Big|\; \sum_{j\in [0,i]}(-1)^j\smatze{i}{j}y_j\con_{t^i} 0\mb{ for all $i\in [0,m-1]$}\Big\} \\
& \tm & \rule[-3mm]{0mm}{0mm} W^{(0)}_m. \\
\end{array}
\;\;$}
\]
\end{Notation}

\begin{Lemma}
\label{LemPT1}
\Absit
\begin{itemize}
\item[(i)] The image $(\Z[\zeta_m] C_m)\omega_m$ of the Wedderburn embedding (\ref{NotSaal0}) is contained in $W^{(1)}_m$,
\[
\Z[\zeta_m] C_m \lraa{\omega_m} W^{(1)}_m \hra W^{(0)}_m.
\]
\item[(ii)] A $\Z[\zeta_m]$-linear basis of $W^{(1)}_m$ is given by 
\[
\left(\xi_{m,i}\right)_{i\in [0,m-1]} := \left(\left((-1)^i t^i\smatze{j}{i}\right)_{j\in [0,m-1]}\right)_{i\in [0,m-1]}.
\]

\item[(iii)] The $(i+1)$st elementary divisor over $\Z[\zeta_m]$ of the embedding $W^{(1)}_m\tm W^{(0)}_m$ is given by $t^i$, $i\in [0,m-1]$. (In particular, $\Z[\zeta_m]$ not being a 
principal ideal domain in general, there exist bases with respect to which the matrix that describes this embedding takes diagonal shape.)

\item[(iv)] The $\Z[\zeta_m]$-linear determinant of the embedding $W^{(1)}_m\tm W^{(0)}_m$ is $t^{m(m-1)/2}$.
\end{itemize}

\rm
Ad (i). We have $\Big((1-X)^i\Big)\left[(\zeta_m^{lj})_{j\in [0,m-1]}\right] =  (1 - \zeta_m^l)^i\con_{t^i} 0$ for all $i\in [0,m-1]$ and all $l\in [0,m-1]$.

Ad (ii). Inverting the matrix $A := \left((-1)^i t^i\smatze{j}{i}\right)_{i\in [0,m-1],\; j\in [0,m-1]}$ arising from the tuple of elements $\left(\xi_{m,j}\right)_{j\in [0,m-1]}$, we 
obtain $A^{-1} := \left((-1)^i t^{-j}\smatze{j}{i}\right)_{i\in [0,m-1],\; j\in [0,m-1]}$. An element $y$ is contained in the $\Z[\zeta_m]$-linear span of our tuple if and only if, $y$
considered as a row vector, $yA^{-1}$ is entrywise contained in $\Z[\zeta_m]$, i.e.\ if and only if $y\in W_m^{(1)}$.
\end{Lemma}

In particular, we have the

\begin{Proposition}
\label{ThW_1}
Let $p$ be a prime. We have a factorization of the Wedderburn embedding (\ref{NotSaal0}) into
\[
\Z[\zeta_p] C_p \lraisoa{\omega_p} W_p^{(1)} \tm W_p^{(0)} = \prod_{j\in \sZ/p}\Z[\zeta_p].
\]

\rm
The factorization follows by (\ref{LemPT1} i). The isomorphism follows by comparison of (\ref{LemPT1} iv) with (\ref{LemIF1}), both yielding the valuation at $t$ of the determinant
of the respective embedding to be $p(p-1)/2$, and zero elsewhere. We remark that the elementary divisors resulting from (\ref{LemPT1} iii) are in accordance with (\ref{PropVand5}).
\end{Proposition}

\begin{quote}
\begin{footnotesize}
\begin{Remark}[coefficient criterion]
\label{RemPT1_1}
The tuple $(y_i)_{i\in [0,p]}\in W_p^{(0)}$ is contained in $W_p^{(1)}$ if and only if, writing $y_i =: \sum_{j\in [0,p-2]} y_{i,j}\zeta^j$, $y_{i,j}\in\Z$,
\[
\sumd{i\in [0,p-1],\; j\in [0,p-2]} (-1)^i\smatze{u}{i} \smatze{j}{v} y_{i,j} \con_p 0
\]
holds for all $0\leq v < u\leq p-1$.

\rm
An element $\sum_{j\in [0,p-2]} x_j\zeta_p^j = \sum_{v\in [0,p-2]}\left(\sum_{j\in [0,p-2]} x_j\smatze{j}{v}\right) (\zeta_p-1)^v$, $x_j\in\Z$, vanishes modulo $t^u$ 
for $u\in [0,p-1]$ if and only if $\sum_{j\in [0,p-2]} x_j\smatze{j}{v} \con_p 0$ for all $v\in [0,u-1]$. It remains to plug in 
$x_j = \sum_{i\in [0,p-1]} (-1)^i \smatze{u}{i} y_{i,j}$. 
\end{Remark}
\end{footnotesize}
\end{quote}

\begin{Lemma}
\label{LemPT2}
Given $0\leq j\leq i$, we have
\[
\xi_{m,j}\xi_{m,i} = \sum_{k\in [0,j]} \smatze{j}{k}\smatze{i+k}{j}(-1)^{j-k}t^{j-k}\xi_{m,i+k},
\]
where we let $\xi_{m,l} := 0$ for $l \geq m$. In particular, $W^{(1)}_m$ is a subring of $W^{(0)}_m$.

\rm
We need to see that for $l\in [0,m-1]$
\[
(-1)^jt^j\smatze{l}{j}\cdot (-1)^i t^i \smatze{l}{i} = \sumd{k\in [0,j]}\smatze{j}{k}\smatze{i+k}{j}(-1)^{j-k}t^{j-k}\cdot (-1)^{i+k}t^{i+k}\smatze{l}{i+k}.
\]
We may assume $i\leq l$ and reformulate to
\[
\smatze{l}{i} = \sumd{k\geq 0} \smatze{j}{k}\smatze{l-j}{l-i-k},
\]
which now follows from a comparison of coefficients in $(1+T)^l = (1+T)^j (1+T)^{l-j}$ at $T^{l-i}$.
\end{Lemma}

\begin{quote}
\begin{footnotesize}
\begin{Remark}
\label{RemPT1_5}
Let $p$ be a prime, let $n\geq 1$, let $m = p^n$ (in particular, $t = 1 - \zeta_{p^n}$). We dispose of ring automorphisms
\[
\begin{array}{rcl}
\Z[\zeta_{p^n}]C_{p^n} & \lraisoah{\alpha_{p^n}} & \Z[\zeta_{p^n}]C_{p^n} \\
c_{p^n}                & \lra                    & \zeta_{p^n} c_{p^n}, \\
\end{array}
\]
and
\[
\begin{array}{rcl}
W_{p^n}^{(0)}       & \lraisoah{\alpha^{(0)}_{p^n}} & W_{p^n}^{(0)} \\
(y_j)_{j\in\sZ/p^n} & \lra                          & (y_{j+1})_{j\in\sZ/p^n} \\
\end{array}
\]
satisfying $\omega_{p^n}\alpha^{(0)}_{p^n} = \alpha_{p^n}\omega_{p^n}$. Moreover, $\alpha^{(0)}_{p^n}$ restricts to an automorphism $\alpha^{(1)}_{p^n}$ of $W_{p^n}^{(1)}$, with operation  
given by
\[
\xi_{p^n,j}\alpha_{p^n}^{(1)} = \xi_{p^n,j} - t\xi_{p^n,j-1} + (-1)^{p-j}t^{j-(p^n-1)}\smatze{p^n}{j}\xi_{p^n,p^n-1}
\]
for $j\in [1,p^n - 1]$, and by $\xi_{p^n,0}\alpha_{p^n}^{(1)} = \xi_{p^n,0} (= 1_{W_{p^n}^{(1)}})$.

\rm
For $j\in [1,p^n-1]$ and $i\in [0,p^n - 1]$, we have
\[
\left((-1)^j t^j\smatze{i}{j}\right) - t\left((-1)^{j-1}t^{j-1}\smatze{i}{j-1}\right) = (-1)^j t^j\smatze{i+1}{j},
\]
which we compare for $i = p^n - 1$ with 
\[
\xi_{p^n,p^n-1} = (0,\dots,0,(-1)^{p^n-1}t^{p^n-1}),
\]
whence the formula describing the operation of $\alpha_{p^n}^{(1)}$. The valuation at $t$ of the third coefficient therein amounts to $j + (n - v_p(j))(p-1)p^{n-1} - (p^n - 1)$. In
case $n - v_p(j)\geq 2$, we obtain
\[
j + (n - v_p(j))(p-1)p^{n-1} \geq 2(p-1)p^{n-1} \geq p^n - 1.
\]
If $n - v_p(j) = 1$, we obtain
\[
j + (n - v_p(j))(p-1)p^{n-1} \geq p^{n-1} + (p-1)p^{n-1} \geq p^n - 1.
\]
Hence $\alpha_{p^n}^{(0)}$ restricts to an automorphism $\alpha_{p^n}^{(1)}$.
\end{Remark}

\begin{Lemma}
\label{LemPT3}
The factorization of the Wedderburn embedding maps
\[
\begin{array}{rcl}
\Z[\zeta_m] C_m & \lraa{\omega_m} & W^{(1)}_m \\
c_m^i           & \lra          & \sum_{k\in [0,m-1]} \left(\frac{1 - \zeta_m^i}{t}\right)^k \xi_{m,k} \\
\end{array}
\]
for $i\in [0,m-1]$, where we let $0^0 = 1$.

\rm
In fact,
\[
\begin{array}{rcl}
\sum_{k\in [0,m-1]} \left(\frac{1 - \zeta_m^i}{t}\right)^k \xi_{m,k} 
& = & \left(\sum_{k\in [0,m-1]} \left(\frac{1 - \zeta_m^i}{t}\right)^k (-1)^k t^k \smatze{j}{k}\right)_{j\in [0,m-1]} \\
& = & \left((1 - (1 - \zeta_m^i))^j\right)_{j\in [0,m-1]} \\
& = & c_m^i\omega_m,
\end{array}
\]
which remains true for $i = 0$.
\end{Lemma}
\end{footnotesize}
\end{quote}
\subsubsection{Second order Pascal ties}
\label{SubSecPasWedSecond}

\begin{footnotesize}
\begin{scriptsize}
\begin{quote}
In this appendix to subsection \ref{SubSecPasWedFirst} we shall indicate a method of how to continue the approach via Pascal ties.
\end{quote}
\end{scriptsize}

Let $n\geq 2$, let $p$ be a prime and let $t = 1 - \zeta_{p^n}$.

\begin{Lemma}
\label{LemPT4}
Let $s\geq 1$. For any $i\geq 1$, we have in $\Q[X]$
\[
\sumd{k\in [1,s]} \frac{1}{k}\left((1-(1-X)^i)^k - i X^k\right) \in (X^{s+1}).
\]

\rm
This expression being contained in $(X)$, it suffices to prove that its derivative is contained in $(X^s)$. In fact,
\[
\begin{array}{cl}
  & \frac{d}{dX}\sumd{k\in [1,s]} \frac{1}{k}\left((1-(1-X)^i)^k - i X^k\right) \\
= & \sumd{k\in [1,s]} \left((1-(1-X)^i)^{k-1}i(1-X)^{i-1} - i X^{k-1}\right) \\
= & i \left(\fracd{1 - (1-(1-X)^i)^s}{1-(1-(1-X)^i)} (1-X)^{i-1} - \fracd{1 - X^s}{1 - X}\right) \\
= & \fracd{i}{1-X} \left(- (1-(1-X)^i)^s + X^s\right). \\
\end{array}
\]
\end{Lemma}

\begin{Lemma}
\label{LemPT5}
For this lemma, we allow $n\geq 1$. Let
\[
f_{p^n}(X) := \left(\sum_{k\in [1,p-1]} \frac{t^{k-1}}{k}\right) X^p - \sum_{k\in [1,p-1]} \frac{t^{k-1}}{k} X^k \in \Z_{(p)}[\zeta_{p^n}][X].
\]
For any $i\in [0,p^n-1]$, we have
\[
f_{p^n}\left(\frac{1 - \zeta_{p^n}^i}{t}\right) \con_{t^{p-1}} 0.
\]
We also write shorthand
\[
\gamma := \sum_{k\in [1,p-1]} \frac{t^{k-1}}{k}.
\]

\rm
In fact,
\[
\begin{array}{ll}
  & \left(\sum_{k\in [1,p-1]} \frac{t^{k-1}}{k}\right) \left(\frac{1 - (1-t)^i}{t}\right)^p
- \sum_{k\in [1,p-1]} \frac{t^{k-1}}{k} \left(\frac{1 - (1-t)^i}{t}\right)^k  \\
= & \frac{1}{t}\sum_{k\in [1,p-1]} \frac{1}{k}\left(t^k\left(\frac{1 - (1-t)^i}{t}\right)^p 
- \left(1 - (1-t)^i\right)^k\right) \\
\con_{p} & \frac{1}{t}\sum_{k\in [1,p-1]} \frac{1}{k}\left(t^k\frac{1 - (1-t^p)^i}{t^p}
- \left(1 - (1-t)^i\right)^k\right) \\
\con_{t^p} & \frac{1}{t}\sum_{k\in [1,p-1]} \frac{1}{k}\left(i t^k - \left(1 - (1-t)^i\right)^k\right) \\
\!\!\auf{\mb{\scr (\ref{LemPT4})}}{\con}_{\!\!\! t^{p-1}} & 0. 
\end{array}
\]
\end{Lemma}

We consider the $\Z[\zeta_{p^n}]$-submodule $W^{(2)}_{p^n}$ of $W^{(1)}_{p^n}$ defined by the {\it second order Pascal ties}
\[
\fbox{$\;\;
\begin{array}{rclcl}
\rule[7mm]{0mm}{0mm} W^{(2)}_{p^n} 
& :=  & \Big\{ \sum_{k\in [0,p^n-1]} z_k\xi_{p^n,k} & \Big| & z_k\in\Z[\zeta_{p^n}],\;\;\left(X^i f_{p^n}^j(X)\right)\left[(z_k)_{k\in [0,p^n-1]}\right]\con_{t^{j(p-1)}} 0 \\
&     &                                       &       & \mb{ for all $i\in [0,p-1]$ and all $j\in [0,p^{n-1}-1]$ } \Big\} \\
& \tm & \rule[-3mm]{0mm}{0mm}W^{(1)}_{p^n}, \\
\end{array}
\;\;$}
\]
where the congruence is to be read in $\Z_{(p)}[\zeta_{p^n}]$.

\begin{Lemma}
\label{LemPT7}
\Absit
\begin{itemize}
\item[(i)] The image $(\Z[\zeta_{p^n}] C_{p^n})\omega_{p^n}$ of the Wedderburn embedding (\ref{NotSaal0}) is contained in $W^{(2)}_{p^n}$,
\[
\Z[\zeta_{p^n}] C_{p^n} \lraa{\omega_{p^n}} W^{(2)}_{p^n} \hra W^{(1)}_{p^n} \hra W^{(0)}_{p^n}.
\]

\item[(ii)] The $(i + jp + 1)$th elementary divisor over $\Z_{(p)}[\zeta_{p^n}]$ of the embedding $(W^{(2)}_{p^n})_{(p)}\tm (W^{(1)}_{p^n})_{(p)}$, $i\in [0,p-1]$, $j\in [0,p^{n-1}-1]$, is given 
by $t^{j(p-1)}$. Moreover, the $(i + jp + 1)$th elementary divisor of the embedding $(W^{(2)}_{p^n})_{(p)}\tm (W^{(0)}_{p^n})_{(p)}$ is given by $t^{(i + jp) + j(p-1)}$.

\item[(iii)] The valuation of the $\Z[\zeta_{p^n}]$-linear determinant of the embedding $W^{(2)}_{p^n}\tm W^{(1)}_{p^n}$ is given by \linebreak
$p^n(p^{n-1}-1)(p-1)/2$ at $t$ and by zero elsewhere.
\end{itemize}

\rm
Ad (i). We have to take care of the coefficients calculated in (\ref{LemPT3}). For $i\in [0,p-1]$, $j\in [0,p^{n-1}-1]$ and $h\in [0,p^n-1]$ we obtain 
\[
\begin{array}{rll}
\left(X^i f_{p^n}^j(X)\right)\left[\left(\left(\frac{1 - \zeta_{p^n}^h}{t}\right)^{\! k}\right)_{k\in [0,p^n-1]}\right]
& = & \left(\frac{1 - \zeta_{p^n}^h}{t}\right)^{\! i} f_{p^n}^j\left(\frac{1 - \zeta_{p^n}^h}{t}\right) \\
& \!\!\!\auf{\mb{\scr (\ref{LemPT5})}}{\con}_{t^{j(p-1)}} & 0.\\
\end{array}
\]
\end{Lemma}

\begin{quote}
\begin{scriptsize}
\begin{Question}
\label{RemPT8}\rm
We do not know whether $W^{(2)}_{p^n}$ is a subring of $W^{(1)}_{p^n}$.
\end{Question}
\end{scriptsize}
\end{quote}

Specializing to $n = 2$, we obtain the

\begin{Proposition}
\label{ThW_2}
We have a factorization of the Wedderburn embedding into
\[
\Z[\zeta_{p^2}] C_{p^2} \lraisoa{\omega_{p^2}} W_{p^2}^{(2)} \tm W_{p^2}^{(1)} \tm W_{p^2}^{(0)} = \prod_{j\in \sZ/p^2}\Z[\zeta_{p^2}].
\]

\rm
The factorization follows by (\ref{LemPT7} i). The isomorphism follows by comparison of (\ref{LemPT7} iii) and (\ref{LemPT1} iv) with (\ref{LemIF1}), both yielding the valuation at $t$ of the 
determinant of the respective embedding to be $p^3(p-1)$, and zero elsewhere. We remark that in this case, the elementary divisors resulting from (\ref{LemPT7} ii) are in accordance with 
(\ref{PropVand5}).
\end{Proposition}

\begin{quote}
\begin{scriptsize}
\begin{Example}
\label{ExPT6}\rm
Consider the case $p = 3$, $n = 2$, thus $t = 1-\zeta_9$, $\gamma = 1 + t/2$. An element $\sum_{k\in [0,8]} z_k\xi_{9,k}$ is contained in $W^{(2)}_9$ if and only if, $z = (z_k)_{k\in [0,8]}$ 
considered as a row vector, it multiplies with
\[
\setlength{\arraycolsep}{2pt}
\renewcommand{\arraystretch}{1.0} 
\left[
\begin{array}{ccc|ccc|ccc}
\scm \; 1\; &\scm 0  &\scm 0    &\scm 0      &\scm 0                &\scm 0                &\scm 0                     &\scm 0                 &\scm 0 \\
\scm 0      &\scm 1  &\scm 0    &\scm -1     &\scm 0                &\scm 0                &\scm 0                     &\scm 0                 &\scm 0 \\
\scm 0      &\scm 0  &\scm 1    &\scm -t/2   &\scm -1               &\scm 0                &\scm 1                     &\scm 0                 &\scm 0 \\\hline
\scm 0      &\scm 0  &\scm 0    &\scm \gamma &\scm -t/2             &\scm -1               &\scm t                     &\scm 1                 &\scm 0 \\
\scm 0      &\scm 0  &\scm 0    &\scm 0      &\scm \gamma           &\scm -t/2             &\scm -2\gamma + t^2/4      &\scm t                 &\scm 1 \\
\scm 0      &\scm 0  &\scm 0    &\scm 0      &\scm 0                &\scm \gamma           &\scm -t\gamma              &\scm -2\gamma + t^2/4  &\scm t \\\hline
\scm 0      &\scm 0  &\scm 0    &\scm 0      &\scm 0                &\scm 0                &\scm \gamma^2              &\scm -t\gamma          &\scm -2\gamma + t^2/4 \\
\scm 0      &\scm 0  &\scm 0    &\scm 0      &\scm 0                &\scm 0                &\scm 0                     &\scm \gamma^2          &\scm -t\gamma \\
\scm 0      &\scm 0  &\scm 0    &\scm 0      &\scm 0                &\scm 0                &\scm 0                     &\scm 0                 &\scm \gamma^2 \\
\end{array}
\weiter
\right]
\left[
\begin{array}{ccc|ccc|ccc}
\scm t^0 &\scm 0   &\scm 0   &\scm 0      &\scm 0      &\scm 0      &\scm 0      &\scm 0      &\scm 0 \\
\scm 0   &\scm t^0 &\scm 0   &\scm 0      &\scm 0      &\scm 0      &\scm 0      &\scm 0      &\scm 0 \\
\scm 0   &\scm 0   &\scm t^0 &\scm 0      &\scm 0      &\scm 0      &\scm 0      &\scm 0      &\scm 0 \\\hline
\scm 0   &\scm 0   &\scm 0   &\scm t^{-2} &\scm 0      &\scm 0      &\scm 0      &\scm 0      &\scm 0 \\
\scm 0   &\scm 0   &\scm 0   &\scm 0      &\scm t^{-2} &\scm 0      &\scm 0      &\scm 0      &\scm 0 \\
\scm 0   &\scm 0   &\scm 0   &\scm 0      &\scm 0      &\scm t^{-2} &\scm 0      &\scm 0      &\scm 0 \\\hline
\scm 0   &\scm 0   &\scm 0   &\scm 0      &\scm 0      &\scm 0      &\scm t^{-4} &\scm 0      &\scm 0 \\
\scm 0   &\scm 0   &\scm 0   &\scm 0      &\scm 0      &\scm 0      &\scm 0      &\scm t^{-4} &\scm 0 \\
\scm 0   &\scm 0   &\scm 0   &\scm 0      &\scm 0      &\scm 0      &\scm 0      &\scm 0      &\scm t^{-4} \\
\end{array}
\right]
\weiter 
\]
to a vector entrywise contained in $\Z_{(3)}[\zeta_9]$.
\end{Example}
\end{scriptsize}
\end{quote}
\end{footnotesize}         
\subsection{Comparison of methods in an example}

\begin{Example}
\label{ExComp1}\rm

Consider the case $m = 5$, in which the cyclic Wedderburn embedding becomes
\[
\begin{array}{rcl}
\Z[\zeta_5]C_5 & \hraa{\omega_5} & \prod_{j\in [0,4]} \Z[\zeta_5] \\
c_5            & \lramaps        & (\zeta_5^j)_{j\in [0,4]}. \\
\end{array}
\]
Let $t = 1-\zeta_5$.

\begin{itemize}
\item[(i)] The $q$-Pascal method (\ref{ThSaal5} iii) yields the $\Z[\zeta_5]$-linear basis of $(\Z[\zeta_5]C_5)\omega_5$ given by the rows of 
\[
\setlength{\arraycolsep}{4pt}
\hspace*{-8mm}
\left[
\begin{array}{ccccc}
5                                     & 0                                      & 0                                         & 0                          & 0 \\
4 + 3\zeta_5 + 2\zeta_5^2 + \zeta_5^3 & -4 - 3\zeta_5 - 2\zeta_5^2 - \zeta_5^3 & 0                                         & 0                          & 0 \\
2 + \zeta_5 + 2\zeta_5^2              & -1 - \zeta_5 + 2\zeta_5^3              & 2 \zeta_5 + \zeta_5^2 + 2 \zeta_5^3       & 0                          & 0 \\
2 + \zeta_5 + \zeta_5^2 + \zeta_5^3   & -2 - \zeta_5^2 - 2 \zeta_5^3           & 3 + 3 \zeta_5 + 2 \zeta_5^2 + 2 \zeta_5^3 & \zeta_5^2 - \zeta_5^3      & 0 \\
1                                     & -1 + \zeta_5 + \zeta_5^3               & - 2 - \zeta_5 - 2 \zeta_5^2 + \zeta_5^3   & - 3 \zeta_5^2 - \zeta_5^3  & \zeta_5 \\
\end{array}
\right].
\]

\item[(ii)] The general Vandermonde method (\ref{LemDiag6}) yields the $\Z[\zeta_5]$-linear basis of $(\Z[\zeta_5]C_5)\omega_5$ given by the rows of 
\[
\hspace*{-8mm}
\setlength{\arraycolsep}{1pt}
\left[
\begin{array}{ccccc}
1 & 1 & 1 & 1 & 1 \\
0 & (\zeta_5-1) & (\zeta_5^2 - 1) & (\zeta_5^3 - 1) & (\zeta_5^4 - 1) \\
0 & 0 & (\zeta_5^2 - 1)(\zeta_5^2 - \zeta_5) & (\zeta_5^3 - 1)(\zeta_5^3 - \zeta_5) & (\zeta_5^4 - 1)(\zeta_5^4 - \zeta_5) \\
0 & 0 & 0 & (\zeta_5^3 - 1)(\zeta_5^3 - \zeta_5)(\zeta_5^3 - \zeta_5^2) & (\zeta_5^4 - 1)(\zeta_5^4 - \zeta_5)(\zeta_5^4 - \zeta_5^2) \\
0 & 0 & 0 & 0 & (\zeta_5^4 - 1)(\zeta_5^4 - \zeta_5)(\zeta_5^4 - \zeta_5^2)(\zeta_5^4 - \zeta_5^3) \\
\end{array}
\right].
\weiter
\]

\item[(iii)] The Pascal method (\ref{LemPT1} ii, \ref{ThW_1}) yields the $\Z[\zeta_5]$-linear basis of $(\Z[\zeta_5]C_5)\omega_5$ given by the rows of 
\[
\left[
\hspace*{-10mm}
\begin{array}{lrrrrr}
\xi_{5,0}\;\; &  1 &     1 &      1 &       1 &      1 \\
\xi_{5,1}     &  0 &    -t &    -2t &     -3t &    -4t \\
\xi_{5,2}     &  0 &     0 &    t^2 &    3t^2 &     6t^2 \\
\xi_{5,3}     &  0 &     0 &      0 &    -t^3 &    -4t^3 \\
\xi_{5,4}     &  0 &     0 &      0 &       0 &      t^4 \\
\end{array}
\right].
\]
But this method works only for $m$ prime (and for $m = p^2$, but less simply).
\end{itemize}
\end{Example}

\section{The radical series of $(W^{(1)}_{p^n})_{(p)}$}

\begin{quote}
\begin{footnotesize}
Let $T$ be a discrete valuation ring with maximal ideal generated by $t$ and residue field $k = T/tT$. Let $\Lambda$ be a subalgebra over $T$ of a direct product of $m$ copies of $T$ such that its 
embedding into this product has torsion cokernel. The dimension over $k$ of the radical layers $\rfk^i\Lambda/\rfk^{i+1}\Lambda$ stabilizes for large $i$ at $m$ [K\"u 99, E.2.3]. On the other 
hand, the surjection
\[
\rfk\Lambda/\rfk^2\Lambda \ts_k \rfk^i\Lambda/\rfk^{i+1}\Lambda \;\lra\;  \rfk^{i+1}\Lambda/\rfk^{i+2}\Lambda
\]
induced by multiplication yields the bound
\[
(\dim_k \rfk\Lambda/\rfk^2\Lambda)(\dim_k \rfk^i\Lambda/\rfk^{i+1}\Lambda) \geq \dim_k(\rfk^{i+1}\Lambda/\rfk^{i+2}\Lambda).
\]
The question is the behaviour of this sequence of dimensions $\dim_k(\rfk^i\Lambda/\rfk^{i+1}\Lambda)$. The ring $W^{(1)}_{p^n}$ introduced in (\ref{NotPT0}), localized 
at $(p)$, yields some example material.
\end{footnotesize}
\end{quote}

\begin{Notation}
\label{NotRS0}\rm
Let $p$ be a prime, let $n\geq 1$ and write $t = 1 - \zeta_{p^n}$. For $i\in\Z$, denote $\ul{i} := \max(0,i)$. Given an integer $i\geq 0$, we write it as $i = \sum_{j\geq 0} a_j p^j$ with
$a_j\in [0,p-1]$ and denote its $p$-adic Quersumme by $q_p(i) := \sum_{j\geq 0} a_j$. 

For $a,b,m\geq 0$, we let
\[
\left(\sum_{i\in [0,m]} T^i\right)^{\!\! a} =: \sum_{b\geq 0} \smatze{a}{b}_m T^b \in \Z[T].
\]
That is, $\smatze{a}{b}_m$ denotes the number of of choices of $b$ balls out of $a$ different balls with each ball chosen at most $m$ times, disregarding the order in which we choose.

In particular, $\smatze{a}{b}_1 = \smatze{a}{b}$. The sequence $(\smatze{a}{b}_m)_{m\geq 0}$ becomes stationary for $m\geq b$ at the value $\smatze{a+b-1}{b}$. We have 
the recursion formula
\[
\smatze{a+1}{b}_m = \sum_{j\in [0,m]} \smatze{a}{b-j}_m,
\]
where we let $\smatze{a}{c}_m := 0$ for $c < 0$. Furthermore, plugging in $T = 1$ in the definition yields $\sum_{b\in [0,am]} \smatze{a}{b}_m = (m+1)^a$.
\end{Notation}

\begin{Proposition}
\label{PropRS1}
\Absit
\begin{itemize}
\item[(i)] The $\Z_{(p)}[\zeta_{p^n}]$-order $(W_{p^n}^{(1)})_{(p)}$ is a local ring, with simple module $\F_p$ acted upon identically by $\xi_{p^n,0}$ and trivially by $\xi_{p^n,i}$ for 
$i\in [1,p-1]$ (cf.\ \ref{LemPT1} ii). I.e.\ 
\[
\rfk (W_{p^n}^{(1)})_{(p)} = (W_{p^n}^{(1)})_{(p)}\cap t (W_{p^n}^{(0)})_{(p)}
\]

\item[(ii)] For $i\geq 0$, the ideal $\rfk^i (W_{p^n}^{(1)})_{(p)}$ is $\Z_{(p)}[\zeta_{p^n}]$-linearly generated by the tuple
\[
\left(t^{\ul{i-q_p(j)}}\xi_{p^n,j}\right)_{j\in [0,p^n-1]}.
\]
In particular, for $i\geq n(p-1)$ we have $\rfk^{i+1} (W_{p^n}^{(1)})_{(p)} = t\rfk^i (W_{p^n}^{(1)})_{(p)}$ (cf.\ [K\"u 99, E.2.3]).

\item[(iii)] On the quantitative side, we obtain for $i\geq 0$
\[
l_{p^n,i} := \dim_{\sF_p} \rfk^i (W_{p^n}^{(1)})_{(p)}/\rfk^{i+1} (W_{p^n}^{(1)})_{(p)} = \sum_{j\in [0,i]} \smatze{n}{j}_{p-1}.
\]
\end{itemize}

\rm
Ad (i). On the one hand, the quotient of the respective $\Z_{(p)}[\zeta_{p^n}]$-order modulo the ideal claimed to be its radical is a simple module. On the other hand, this ideal is 
nilpotent modulo $t$.

Ad (ii). We abbreviate $\rfk^i := \rfk^i (W_{p^n}^{(1)})_{(p)}$ for $i\geq 0$ and perform an induction on $i$, the case of $i = 1$ being true by (i), and a basis of $\rfk\Lambda$ being given 
by $(t\xi_{p^n,0},\xi_{p^n,1},\dots,\xi_{p^n,p^n-1})$. Note that for $0\leq b\leq a$, the binomial coefficient $\smatze{a}{b}$ has valuation $p^{n-1}(q_p(a-b) + q_p(b) - q_p(a))$ at $t$.

{\it We claim that $\rfk\cdot\rfk^{i-1}$ is contained in the span of the given tuple.}

The set $t\xi_{p^n,0}\cdot\rfk^{i-1}$ is contained in this span. Moreover, by (\ref{LemPT2}), we obtain in case $1\leq j'$, $0\leq j\leq j'$, 
$$
\xi_{p^n,j'}(t^{\ul{i-1-q_p(j)}}\xi_{p^n,j}) = t^{\ul{i-1-q_p(j)}} \sum_{s\in [0,j]} \smatze{j}{s}\smatze{j+j'-s}{j} (-1)^st^s \xi_{p^n,j+j'-s}.
\leqno (\ast)
$$
So it suffices to show the inequality
\[
\begin{array}{rcl}
v_t(\smatze{j+j'-s}{j}) + \ul{i-1-q_p(j)} + s & \geq & \ul{i-q_p(j+j'-s)} \\
                                                & =    & \ul{i - q_p(j) - q_p(j'-s) + p^{1-n}v_t(\smatze{j+j'-s}{j})} \\
\end{array}
\]
for $s\in [0,j]$, whence in turn it suffices to show that
\[
(1-p^{1-n})v_t(\smatze{j+j'-s}{j}) - 1 + s + q_p(j'-s) \geq 0.
\]
But we have $s\geq 1$ or ($s = 0$ and $q_p(j')\geq 1$).

Furthermore, (\ref{LemPT2}) yields
$$
\xi_{p^n,j'}(t^{\ul{i-1-q_p(j)}}\xi_{p^n,j}) = t^{\ul{i-1-q_p(j)}} \sum_{s\in [0,j']} \smatze{j'}{s}\smatze{j+j'-s}{j'} (-1)^st^s \xi_{p^n,j+j'-s}.
\leqno (\ast\ast)
$$
in case $1\leq j'\leq j$. So it suffices to show that
\[
v_t(\smatze{j'}{s}\smatze{j+j'-s}{j'}) + \ul{i-1-q_p(j)} + s \geq \ul{i-q_p(j+j'-s)}
\]
for $s\in [0,j']$, whence in turn, dropping a factor $p^{n-1}$, it suffices to see that
\[
\Big((q_p(j'-s) + q_p(s) - q_p(j')) + (q_p(j-s) + q_p(j') - q_p(j+j'-s))\Big) - q_p(j) + q_p(j+j'-s) + s - 1\geq 0.
\]
For $s\geq 1$, we see that both in case $(\ast)$ and in case $(\ast\ast)$, the corresponding summand is even contained in $t\rfk^{i - 1}$.

{\it We claim that $\rfk\cdot\rfk^{i-1}$ contains the given tuple.}

First, we note that $t^{\ul{i-q_p(0)}}\xi_{p^n,0} = (t\xi_{p^n,0})^i$. Moreover, by the summandwise argument above concerning $(\ast\ast)$, we dispose of the congruence
\[
\xi_{p^n,j'}(t^{\ul{i-1-q_p(j)}}\xi_{p^n,j}) \con_{t\rfk^{i-1}} t^{\ul{i-1-q_p(j)}} \smatze{j+j'}{j'} \xi_{p^n,j+j'}.
\]
in case $1\leq j'\leq j$. In particular, given $j\in [1,p^n-1]$ such that $j[p]\leq j - j[p]$, i.e.\ such that $j$ is not a power of $p$, we obtain
\[
\xi_{p^n,j[p]}(t^{\ul{i-1-q_p(j-j[p])}}\xi_{p^n,j-j[p]}) \con_{t\rfk^{i-1}} t^{\ul{i-q_p(j)}} \smatze{j}{j[p]} \xi_{p^n,j}\; ,
\]
and it remains to note that $v_p(\smatze{j}{j[p]}) = 0$. 

Furthermore, for $m\in [0,n-1]$ we have
\[
\xi_{p^n,p^m}(t^{\ul{i-1-q_p(0)}}\xi_{p^n,0}) = t^{\ul{i-q_p(p^m)}} \xi_{p^n,p^m}.
\]

Ad (iii). Concerning the dimension of the $i$th radical layer, we calculate
\[
\begin{array}{rcl}
\#\{k\in [0,p^n-1]\; |\; q_p(k) = j\} 
& = & \#\{ (a_l)_{l\in [0,n-1]}\in [0,p-1]^n\; |\; \sum_{l\in [0,n-1]} a_l = j\} \\
& = & \smatze{n}{j}_{p-1}. \\
\end{array}
\]
\end{Proposition}

\begin{quote}
\begin{footnotesize}
\begin{Remark}
\label{RemRS2}
\rm
The inequality mentioned in the introduction to this section reads $l_{p^n,1}\cdot l_{p^n,i}\geq l_{p^n,i+1}$ for $i\geq 1$. By (\ref{PropRS1} iii), this translates into the assertion that
\[
(n+1)\left(\sum_{j\in [0,i]} \smatze{n}{j}_{p-1}\right) \geq \sum_{j\in [0,i+1]} \smatze{n}{j}_{p-1}.
\]
For $p$ large, this follows from
\[
\begin{array}{rcl}
(n+1)\left(\sum_{j\in [0,i]} \smatze{n}{j}_{p-1}\right)
& = & (n+1)\left(\sum_{j\in [0,i]} \smatze{n + j - 1}{n - 1}\right) \\
& = & (n+1)\smatze{n+i}{n} \\
& \geq & \frac{n+i+1}{i+1}\cdot\frac{(n+i)!}{n!i!} \\
& = & \smatze{n+i+1}{n} \\
& = & \sum_{j\in [0,i+1]} \smatze{n}{j}_{p-1}.\\
\end{array}
\]
\end{Remark}

\begin{Example}
\label{ExRS3}
The sequence $(l_{81,i})_{i\geq 0}$ is given by
\[
1,\, 5,\, 15,\, 31,\, 50,\, 66,\, 76,\, 80,\, 81,\, 81,\,\dots
\]
\end{Example}
\end{footnotesize}
\end{quote}
\section{The absolute cyclic Wedderburn embedding}
\label{SecAbs}

\begin{quote}
\begin{footnotesize}
Once given the Kervaire-Murthy pullback, a closed formula describing the image of the absolute Wedderburn embedding can be derived in a straightforward manner (\ref{ThKM5}, \ref{PropCP9_2}).
\end{footnotesize}
\end{quote}

\subsection{An inversion formula for $\Q C_{p^n}$}
\label{SubSecInvpn}

Let $p$ be a prime and let $n\geq 1$.

\begin{Lemma}
\label{LemCP9}
The inverse of the absolute Wedderburn isomorphism
\[
\begin{array}{rcl}
\Q C_{p^n} & \lraisoa{\omega_{\ssQ,p^n}} & \prodd{k\in [0,n]} \Q[\zeta_{p^k}] \\
c_{p^n}    & \lramaps                    & (\zeta_{p^k})_{k\in [0,n]} \\
\end{array}
\]
is given by
\[
\begin{array}{rcl}
\prodd{k\in [0,n]} \Q[\zeta_{p^k}] & \lraisoah{\omega_{\ssQ,p^n}^{-1}} & \Q C_{p^n} \\
(\sumd{i\in\sZ} y_{k,i}\zeta_{p^k}^i)_{k\in [0,n]} & \lramaps & 
p^{-n}\sumd{i\in\sZ}\left(\sumd{k\in [0,n]} p^k y_{k,i}\!\!\!\!\sumd{j\in [1,p^{n-k}]} c_{p^n}^{i+jp^k} 
      - \sumd{k\in [1,n]} p^{k-1}y_{k,i}\!\!\!\!\!\sumd{j\in [1,p^{n-k+1}]} c_{p^n}^{i+jp^{k-1}}\right) \\
& = & p^{-n}\sumd{i\in\sZ}c_{p^n}^i \sumd{k\in [0,n]} p^k\sumd{j\in [1,p^{n-k}]}  \left(y_{k,i-jp^k} - y_{k+1,i-jp^k}\right), \\
\end{array}
\]
where $y_{k,i}\in\Q$ and where $y_{n+1,i} := 0$ for $i\in\Z$. We allow `non-reduced' expressions for elements in $\Q(\zeta_{p^k})$, merely requiring $y_{k,i} = 0$ for all but finitely 
many $i\in\Z$. 

We denote the absolute Wedderburn embedding, i.e.\ the restriction of $\omega_{\sQ,p^n}$ to $\Z C_{p^n}$, by
\[
\Z C_{p^n} \lraa{\omega_{\ssZ,p^n}} \prodd{k\in [0,n]} \Z[\zeta_{p^k}].
\]

\rm
By (\ref{LemIF1}), the inversion formula over $\Q(\zeta_{p^n})$ reads 
\[
\begin{array}{rcl}
\Q(\zeta_{p^n}) C_{p^n} & \lraiso  & \prodd{i\in [1,p^n]} \Q(\zeta_{p^n}) \\
c_{p^n}                 & \lramaps & (\zeta_{p^n}^i)_{i\in [1,p^n]} \\ 
p^{-n}\sumd{j\in [1,p^n]}c_{p^n}^j\sumd{i\in [1,p^n]}y_i\zeta_{p^n}^{-ij} & \llamaps & (y_i)_i. \\
\end{array}
\]

We precompose the inverse direction with the direct product of the maps
\[
\begin{array}{rcr}
\Q(\zeta_{p^k}) & \lra     & \prodd{\sigma\in\sGal(\sQ(\zeta_{p^k})/\sQ)} \Q(\zeta_{p^n}) \\
x               & \lramaps & (x\sigma)_\sigma
\end{array}
\]
over $k\in [0,n]$, where the position of the Galois automorphism corresponding to $i\in (\Z/p^k)^\ast$, represented by $i\in [1,p^k]$, is to be identified with the position
$i p^{n-k} \in [1,p^n]$. 

If $l\in [1,n]$, the image of $(\dell_{l,k})_{k\in [0,n]}\in\prod_{k\in [0,n]}\Z[\zeta_{p^k}]$ under this composition is 
\[
\begin{array}{rl}
  & p^{-n}\sumd{j\in [1,p^n]} c_{p^n}^j\sumd{i\in [1,p^n],\; v_p(i)\; =\; n-l} \zeta_{p^n}^{-ij} \\
= & p^{-n}\sumd{j\in [1,p^n]} c_{p^n}^j\sumd{i\in (\sZ/p^l)^\ast} \zeta_{p^n}^{-ijp^{n-l}} \\ 
= & p^{-n}\left(\sumd{j\in [1,p^n],\; v_p(j)\geq l} c_{p^n}^j (p-1)p^{l-1}  
+   \sumd{j\in [1,p^n],\; v_p(j)\; =\; l - 1} c_{p^n}^j (-1) p^{l-1}\right) \\ 
= & p^{-n}\left(p^l\left(\sumd{j\in [1,p^{n-l}]} c_{p^n}^{jp^l}\right)
-   p^{l-1}\left(\sumd{j\in [1,p^{n-l+1}]} c_{p^n}^{jp^{l-1}}\right)\right). \\
\end{array}
\]
If $l = 0$, it is $p^{-n}\sumd{j\in [1,p^n]} c_{p^n}^j$. The inversion formula follows by $\Q C_{p^n}$-linearity of $\omega_{\sQ,p^n}$.
\end{Lemma}

\subsection{An inversion formula for $\Q C_m$}
\label{SubSecInvm}

Let $m\geq 1$. Writing $d|m$ stands for $d\geq 1$ and $m\in (d)$. The letters $p$ and $q$ are reserved to denote prime numbers. We denote by $k[p'] := k/k[p]$ the $p'$-part of an integer $k\geq 1$.
Given a finite family of groups $(G_i)_{i\in [1,l]}$, $l\geq 1$, there is a ring isomorphism
\[
\begin{array}{rcl}
\Z[\prodd{i\in [1,l]}G_i] & \lraiso  & \Tsd{i\in [1,l]} \Z G_i \\
(g_i)_i                   & \lramaps & (g_i)^\ts_i := g_1\ts \cdots \ts g_l. \\
\end{array}
\]
For $p|k$ we denote by $s_{k,p}$ a representative in $\Z$ of the inverse of $k[p']$ in $(\Z/k[p])^\ast$ to obtain
\[
\begin{array}{rcl}
C_k                               & \lraiso  & \prodd{p|k} C_{k[p]} \\
c_k                               & \lramaps & (c_{k[p]})_{p|k} \\
c_k^{\sum_{p|k} j_p s_{k,p}k[p']} & \llamaps & (c_{k[p]}^{j_p})_{p|k}, \\
\end{array}
\]
and
\[
\begin{array}{rcl}
\Z[\zeta_k]                            & \llaiso  & \Tsd{p|k} \Z[\zeta_{k[p]}] \\
\zeta_k^{\sum_{p|k} s_{k,p} j_p k[p']} & \llamaps & (\zeta_{k[p]}^{j_p})^\ts_{p|k} \\
\zeta_k                                & \lramaps & (\zeta_{k[p]})^\ts_{p|k}. \\
\end{array}
\]
We use these isomorphisms as identifications. 

We consider the absolute Wedderburn isomorphism
\[
\begin{array}{rcl}
\Q C_m & \lraisoa{\omega_{\ssQ,m}} & \prodd{d|m} \Q(\zeta_d) \\
c_m    & \lramaps                  & (\zeta_d)_{d|m} \\
\end{array}
\]
and its restriction, the absolute Wedderburn embedding 
\[
\Z C_m \lraa{\omega_{\ssZ,m}} \prodd{d|m} \Z[\zeta_d].
\]

\begin{Remark}
\label{RemarkCP9_0}
The embedding $\omega_{\sZ,m}$ identifies with $\Ts_{p|m} \omega_{\sZ,m[p]}$ in the sense that the identifications yield an isomorphism of embeddings
\begin{center}
\begin{picture}(400,250)
\put(-170, 200){$\Z C_m$}
\put( -70, 210){\vector(1,0){450}}
\put(  50, 225){$\scm\omega_{\ssZ,m}$}
\put( 400, 200){$\prodd{d|m}\Z[\zeta_d]$}
\put(-125, 180){\vector(0,-1){130}}
\put( 450,  50){\vector(0,1){90}}
\put(-230,   0){$\Ts_{p|m}\Z C_{m[p]}$}
\put(  10,  10){\vector(1,0){270}}
\put(  50,  30){$\scm\Ts_{p|m}\omega_{\ssZ,m[p]}$}
\put( 300,   0){$\Ts_{p|m}\prodd{e|m[p]}\Z[\zeta_e]$.}
\end{picture}
\end{center}

\rm
On both ways, $c_m$ is mappped to $(\zeta_d)_{d|m}$.
\end{Remark}

We factorize $\omega_{\sQ,m}^{-1}$ along identifications $u$ and $w$ into
\[
\prodd{d|m} \Q(\zeta_d) \lraisoa{u} \Tsd{p|m} \left(\prodd{d|m[p]} \Q[\zeta_d]\right) \lraisoa{v} \Ts_{p|m} \Q C_{m[p]} \lraisoa{w} \Q C_m,
\]
$v$ being the tensor product of the inverses of the absolute Wedderburn isomorphism belonging to the respective prime part (cf.\ \ref{RemarkCP9_0}). 

For $e|m$ and a prime $p|m$, we abbreviate the element
\[
\Q C_m \ni f_{e,p}(c_m) := 
\left\{
\begin{array}{ll}
\frac{1}{m[p]}\sumd{i\in\sZ/m[p]} c_m^i & \mb{for } v_p(e) = 0 \\
\frac{e[p]}{m[p]}\left(\sumd{i\in\sZ/(\frac{m[p]}{e[p]})} c_m^{ie[p]} - \frac{1}{p}\sumd{i\in\sZ/(p\frac{m[p]}{e[p]})} c_m^{ie[p]/p}\right) 
                                               & \mb{for } v_p(e) \geq 1. \\
\end{array}
\right.
\]
Suppose given an element $a_d(c_m) = \sumd{i\in \sZ/m} a_{d,i} c_m^i\in\Q C_m$ for each $d|m$, representing $a_d(\zeta_d) = \sum_{i\in\sZ/m} a_{d,i}\zeta_d^i$. Using $\Q C_m$-linearity, we obtain
\[
\begin{array}{rcl}
(a_d(\zeta_d))_{d|m} uvw
& = & \sum_{e|m} a_e(c_m)\cdot (\dell_{e,d})_{d|m} uvw \\
& = & \sum_{e|m} a_e(c_m)\cdot\left( (\dell_{e[p],d})_{d|m[p]} \right)^\ts_{p|m} vw \\
& \aufgl{(\ref{LemCP9})} & \sum_{e|m} a_e(c_m)\cdot\left( f_{e,p}(c_{m[p]}) \right)^\ts_{p|m} w \\
& = & \sum_{e|m} a_e(c_m) \prod_{p|m} f_{e,p}(c_m^{s_{m,p} m[p']}). \\
\end{array} 
\]
We write 
\[
f_{e,p}(c_m) =: \sumd{j\in\sZ/m[p]} f_{e,p,j} c_m^{j}
\]
and continue to calculate
\[
\begin{array}{rcl}
\sum_{e|m} a_e(c_m) \prod_{p|m} f_{e,p}(c_m^{s_{m,p} m[p']})
& = & \sum_{e|m}\sum_{i\in\sZ/m} a_{e,i} c_m^i \prod_{p|m} \sum_{j\in\sZ/m[p]}f_{e,p,j} c_m^{j s_{m,p} m[p']} \\
& = & \sum_{e|m}\sum_{i\in\sZ/m} a_{e,i} c_m^i \sum_{k\in\sZ/m} c_m^k\prod_{p|m} f_{e,p,k} \\
& = & \sum_{l\in\sZ/m} c_m^l\left[\sum_{e|m}\sum_{k\in\sZ/m} a_{e,l-k}\prod_{p|m} f_{e,p,k}\right] \\
\end{array}
\]
If $A$ is a condition which might be true or not, we let the expression $\{\mb{if } A \}$ take the value $1$ if $A$ holds, and $0$ if
$A$ does not hold. Given $d|m$, we let $d'|m$ be defined by
\[
v_p(d') := \max(v_p(d) - 1,0)
\]
for $p|m$. We rewrite
\[
\text
f_{e,p,k} = \frac{e[p]}{m[p]}\left[\{\mb{if }e[p]|k\} - \frac{1}{p}\{\mb{if } p|e[p]\}\{\mb{if } e[p]|pk\} \right]
\]
to obtain
\[
\begin{array}{rcl}
\prod_{p|m} f_{e,p,k} 
& = & \frac{e}{m} \prod_{p|m}\left[\{\mb{if }e[p]|k\} - \frac{1}{p}\{\mb{if } p|e[p]\}\{\mb{if } e[p]|pk\} \right] \\ 
& = & \{\mb{if } e' | k \} \cdot\frac{e}{m} 
\left[\prodd{p|e,\; e[p] = pk[p]} \left(-\frac{1}{p}\right)\right]\left[\prodd{p|e,\; e[p]|k[p]} \left(1-\frac{1}{p}\right)\right]. \\
\end{array}
\]

\begin{Proposition}[inversion formula]
\label{PropCP10}
The inverse of the absolute Wedderburn isomorphism
\[
\begin{array}{rcl}
\Q C_m & \lraisoa{\omega_{\ssQ,m}} & \prodd{d|m} \Q(\zeta_d) \\
c_m    & \lramaps                  & (\zeta_d)_{d|m} \\
\end{array}
\]
maps $(a_d(\zeta_d))_{d|m}$, written using representing elements $a_d(c_m) = \sumd{i\in \sZ/m} a_{d,i} c_m^i\in\Q C_m$, to 
\[
\frac{1}{m}\sum_{l\in\sZ/m} c_m^l\left[\sum_{d|m}d\sum_{k\in\sZ/(\frac{m}{d'})} a_{d,l-kd'}\left[\prodd{p|d,\; p\!\not\;|\, k} \left(-\frac{1}{p}\right)\right]
\left[\prodd{p|d,\; p|k} \left(1-\frac{1}{p}\right)\right]\right],
\]
where $v_p(d') := \max(v_p(d) - 1,0)$.
\end{Proposition}

\begin{quote}
\begin{footnotesize}
But these inversion formulas (\ref{LemCP9}, \ref{PropCP10}) are not suited for giving convenient ties that describe the image of the absolute Wedderburn embedding, for in neither case they yield
a {\it triangular} system of ties, i.e.\ of congruences of tuple entries.
\end{footnotesize}
\end{quote}

\subsection{The Kervaire-Murthy pullback}
\label{SubSecKerMurPB}

Let $p$ be a prime and let $n\geq 1$.

\begin{Lemma}
\label{LemKM1}
Let $k\in [1,n]$. Writing
\[
\begin{array}{rcl}
f_k(X) & := & -\sumd{i\in [0,p-2]} (p-1-i) X^{ip^{k - 1}} \\
g_k(X) & := & \sumd{i\in [0,p^{n-k}-1]} \left(- (p^{n-k+1}-p(i+1))     X^{ip^k + p^{k - 1}} 
                                                           + (p^{n-k+1}-p(i+1) + 1) X^{ip^k}\right), \\
\end{array}
\]
we obtain
\[
  f_k(X)\cdot \left(\prodd{i\in [0,n]\ohne\{k\}} \Phi_{p^i}(X)\right)
+ g_k(X)\cdot \Phi_{p^k}(X)
= p^{n - k + 1}.
\]

\rm
We expand
\[
\begin{array}{rl}
  & f_k(X)\cdot \left(\prodd{i\in [0,n]\ohne\{k\}} \Phi_{p^i}(X)\right) \\
= & -\left(\sum_{i\in [0,p-2]} (p-1-i) X^{ip^{k - 1}}\right)\left((X^{p^{k-1}} - 1)\sum_{i\in [0,p^{n-k}-1]} X^{ip^k}\right) \\
= & -\left(\sum_{i\in [1,p-1]} (p-i) X^{ip^{k - 1}} - \sum_{i\in [0,p-2]} (p-1-i) X^{ip^{k - 1}}\right)
    \left(\sum_{i\in [0,p^{n-k}-1]} X^{ip^k}\right) \\
= & -\left(-p + \sum_{i\in [0,p-1]} X^{ip^{k - 1}}\right)
    \left(\sum_{i\in [0,p^{n-k}-1]} X^{ip^k}\right) \\
= & p \left(\sum_{i\in [0,p^{n-k}-1]} X^{ip^k}\right) - \left(\sum_{i\in [0,p^{n-k+1}-1]} X^{ip^{k-1}}\right). \\
\end{array} 
\]
Then we calculate
\[
\begin{array}{rl}
  & g_k(X)\cdot \Phi_{p^k}(X) \\
= & \left(\sum_{i\in [0,p^{n-k}-1]} \left(- (p^{n-k+1}-p(i+1))     X^{ip^{k} + p^{k - 1}} 
                                                   + (p^{n-k+1}-p(i+1) + 1) X^{ip^{k}}\right)\right)\cdot \\
  & \cdot\left(\sum_{j\in [0,p-1]} X^{jp^{k-1}}\right)\\
= & -\sum_{i\in [0,p^{n-k}-1]}\sum_{j\in [1,p]}    (p^{n-k+1}-p(i+1))     X^{ip^k + jp^{k-1}} \\
  & +\sum_{i\in [0,p^{n-k}-1]}\sum_{j\in [0,p-1]}  (p^{n-k+1}-p(i+1) + 1) X^{ip^k + jp^{k-1}} \\
= & -\sum_{i\in [0,p^{n-k}-1]} (p^{n-k+1}-p(i+1)) X^{(i+1)p^k} 
    +\sum_{i\in [0,p^{n-k}-1]}\sum_{j\in [1,p-1]}  X^{ip^k + jp^{k-1}} \\
  & +\sum_{i\in [0,p^{n-k}-1]} (p^{n-k+1}-p(i+1) + 1) X^{ip^k} \\
= & -\sum_{i\in [1,p^{n-k}-1]} (p^{n-k+1}-pi) X^{ip^k} 
    +\sum_{i\in [0,p^{n-k+1}-1]} X^{ip^{k-1}}  \\
  & +\sum_{i\in [0,p^{n-k}-1]} (p^{n-k+1}-p(i+1)) X^{ip^k} \\
= & -p\left(\sum_{i\in [1,p^{n-k}-1]} X^{ip^k}\right) 
    + \left(\sum_{i\in [0,p^{n-k+1}-1]} X^{ip^{k-1}}\right)  
    + (p^{n-k+1}-p).  \\
\end{array}
\]
\end{Lemma}

\begin{Lemma}
\label{LemKM2}
Writing
\[
g_0(X) := -\sumd{i\in [0,p^n-2]} (p^n-1-i) X^i,
\]
we obtain
\[
\left(\prodd{i\in [0,n]\ohne\{ 0\}}\Phi_{p^i}(X)\right)
+ g_0(X)\cdot\Phi_{p^0}(X)
= p^n.
\]

\rm
We expand
\[
\begin{array}{rcl}
g_0(X)\cdot\Phi_{p^0}(X) 
& = & - \left(\sum_{i\in [0,p^n-2]} (p^n-1-i) X^i\right)(X-1) \\
& = & - \left(\sum_{i\in [1,p^n-1]} (p^n-i) X^i\right) + \left(\sum_{i\in [0,p^n-1]} (p^n-1-i) X^i\right) \\
& = & - \left(\sum_{i\in [0,p^n-1]} X^i \right) +  p^n.  \\
\end{array}
\]
\end{Lemma}

\begin{Proposition}[{{\sc Kervaire, Murthy} [KM 77, \S 1]}]
\label{PropKM3}
Let $k\geq 1$. There is a pullback diagram of rings
\begin{center}
\begin{picture}(250,350)
\put( -60, 250){$\Z C_{p^k}$}
\put(  50, 260){\vector(1,0){130}}
\put( 200, 250){$\Z C_{p^{k-1}}$}
\put( -20, 230){\vector(0,-1){130}}
\put( -60, 170){$\scm\alpha_k$}
\put( 220, 230){\vector(0,-1){130}}
\put( -65,  50){$\Z[\zeta_{p^k}]$}
\put(  50,  60){\vector(1,0){130}}
\put( 200,  50){$\F_p C_{p^{k-1}}$}
\put(  50, 210){\line(1,0){40}}
\put(  50, 210){\line(0,-1){40}}
\put(-110, 305){$\scm c_{p^k}$}
\put( -60, 310){\vector(1,0){400}}
\put( -60, 300){\line(0,1){20}}
\put( 363, 305){$\scm c_{p^{k-1}}$}
\put( -90, 280){\vector(0,-1){250}}
\put(-100, 280){\line(1,0){20}}
\put( 370, 280){\vector(0,-1){250}}
\put( 360, 280){\line(1,0){20}}
\put(-117,   0){$\scm\zeta_{p^k}$}
\put( -60,  10){\vector(1,0){400}}
\put( -60,   0){\line(0,1){20}}
\put( 363,   0){$\scm c_{p^{k-1}}$.}
\end{picture}
\end{center}

\rm
Given a commutative ring $A$ containing ideals $\afk$ and $\bfk$, there is a pullback diagram
\begin{center}
\begin{picture}(250,250)
\put(-100, 200){$A/\afk\cap\bfk$}
\put(  50, 210){\vector(1,0){130}}
\put( 200, 200){$A/\bfk$}
\put(  20, 180){\vector(0,-1){130}}
\put( 220, 180){\vector(0,-1){130}}
\put( -40,  00){$A/\afk$}
\put(  50,  10){\vector(1,0){130}}
\put( 200,  00){$A/\afk+\bfk$.}
\end{picture}
\end{center}

Hence, letting $A = \Z[X]$, $\afk = (\Phi_{p^k}(X))$, $\bfk = (X^{p^{k-1}}-1)$, we need to show that 
$(\Phi_{p^k}(X))\cap (X^{p^{k-1}}-1) = (X^{p^k}-1)$ and that $(\Phi_{p^k}(X),X^{p^{k-1}}-1) = (p,X^{p^{k-1}}-1)$. The intersection
follows from the analoguous assertion in $\Q[X]$, using $\Z[X]f(X) = \Z[X]\cap\Q[X]f(X)$ for a monic polynomial $f(X)\in\Z[X]$. The sum follows by 
\[
(\Phi_{p^k}(X),X^{p^{k-1}}-1)\aufgl{(\ref{LemKM1})} (p,\Phi_{p^k}(X),X^{p^{k-1}}-1) 
= (p,(X-1)^{p^k-p^{k-1}},(X-1)^{p^{k-1}}).
\]
\end{Proposition}

\begin{Corollary}
\label{CorKM4}
The index of the absolute Wedderburn embedding $\omega_{\sZ,p^n}$ is given by 
\[
\#\prodd{k\in [1,n]}\F_p C_{p^{k-1}} = p^{\frac{p^n-1}{p-1}}. 
\]
\end{Corollary}
\begin{Corollary}
\label{CorCP9_0_5}
For $m\geq 1$, the index of the absolute Wedderburn embedding $\omega_{\sZ,m}$ is given by
\[
\prodd{p|m}\; p^{\frac{m[p]-1}{p-1}\cdot m[p']}.
\]

\rm
This follows from (\ref{CorKM4}) using (\ref{RemarkCP9_0}).
\end{Corollary}

\begin{quote}
\begin{footnotesize}
\begin{Remark}
\label{RemCP9_0_6}
\rm
For $k\geq 1$, we let $\Delta_k$ denote the absolute value of the discriminant of $\Z[\zeta_k]$ over $\Z$. We take from [K\"u 99, S 1.1.2] that the index of $\omega_{\sZ,m}$ is given by
\[
\sqrt{\frac{m^m}{\prod_{d|m}\Delta_d}}.
\]
A comparison with (\ref{CorCP9_0_5}) allows to re-calculate $\Delta_m$. First we remark that the inverse of the Dedekind isomorphism
\[
\begin{array}{rclcl}
\Q(\zeta_m) & \ts & \Q(\zeta_m) & \lraisoa{\delta_m} & \prod_{j\in(\sZ/m)^\ast}\Q(\zeta_m) \\
\zeta_m^k   & \ts &\zeta_m^l    & \lramaps           & (\zeta_m^k\zeta_m^{jl})_{j\in(\sZ/m)^\ast}, \\
\end{array}
\]
where $k,l\in [0,\phi(m)-1]$, is given by
\[
\begin{array}{rcl}
\prod_{j\in (\sZ/m)^\ast}\Q(\zeta_m) & \lraa{\delta_m^{-1}} & \Q(\zeta_m)\ts\Q(\zeta_m) \\
(y_j)_{j\in (\sZ/m)^\ast}            & \lramaps             & m^{-1} \sumd{i\in\sZ/m}\;\;\sumd{j\in(\sZ/m)^\ast} y_j\zeta_m^{-ij}\ts\zeta_m^i. \\
\end{array}
\]
In particular, its restriction to the Dedekind embedding
\[
\Z[\zeta_m] \ts \Z[\zeta_m] \hraa{\delta_m} \prod_{j\in(\sZ/m)^\ast}\Z[\zeta_m]
\]
has a cokernel annihilated by $m$. In particular, the prime divisors of $\Delta_m$ divide $m$. We infer that for $p|m$
\[
\sumd{d|m,\; v_p(d)\geq 1} v_p(\Delta_d) = \sumd{d|m}v_p(\Delta_d) 
\aufgl{(\ref{CorCP9_0_5})} m[p']\left(v_p(m)m[p] - 2\cdot\frac{m[p]-1}{p-1}\right).
\]
Given $d\geq 1$, $p|d$, we denote
\[
u(p,d) := \frac{d[p]}{p}\left(\rule{0mm}{4mm}v_p(d)(p-1)-1\right)\cdot\phi(d[p'])
\]
and obtain
\[
\begin{array}{rcl}
\sumd{d|m,\; v_p(d)\geq 1} u(p,d)
& = & \sumd{d|m,\; v_p(d)\geq 1}\frac{d[p]}{p}\left(\rule{0mm}{4mm}v_p(d)(p-1)-1\right)\cdot\phi(d[p']) \\
& = & \sumd{i\in [1,v_p(m)]}p^{i-1}\left(\rule{0mm}{4mm}i (p-1)-1\right)\cdot\sumd{d|m[p']}\phi(d) \\
& = & \left(v_p(m)m[p] - 2\frac{m[p]-1}{p-1}\right)\cdot m[p']\; , \\
\end{array}
\]
whence by induction on $m$, we get $v_p(\Delta_m) = u(p,m)$, i.e.\ the absolute value of the \linebreak discriminant is
{\ncr
\[
\Delta_m = \prodd{p|m}\; p^{\frac{m[p]}{p}\left(\rule{0mm}{4mm}v_p(m)(p-1)-1\right)\phi(m[p'])} \; .
\]
}
\end{Remark}
\end{footnotesize}
\end{quote}

\subsection{Kervaire-Murthy ties for $\Z C_{p^n}$}
\label{SubSecKerMurTies}

Let $p$ be a prime and let $n\geq 0$.
Given $j\in\Z$ and $k\geq 1$, we let $[j]_k\in [0,k-1]$ be defined by $[j]_k\con_k j$. For $a,b,c\in\Z$ and $a\leq b$, we let $\chi_{[a,b]}(c) := 1$ if $c\in [a,b]$ and $\chi_{[a,b]}(c) := 0$ 
if $c\not\in [a,b]$.

Let $\lfk^\mb{\scr\rm f}(\Z)$ be the abelian group consisting of sequences $x = (x_j)_{j\in\sZ}$ with entries $x_j\in\Z$ such that the support 
\[
\ul{x} := \{ i\in\Z \; |\; x_i\neq 0\}\tm\Z
\]
of positions carrying nonzero entries is finite. For $i\geq 0$, $m\geq 0$ and $s\geq 0$ we shall define a $\Z$-linear operator $T^{m,s}_i: \lfk^\mb{\scr\rm f}(\Z) \lra \lfk^\mb{\scr\rm f}(\Z)$ 
(which we write on the {\it left}).

Suppose given $j\in\Z$. If $i\geq m$, we let
\[
\Big(T^{m,s}_i x\Big)_j := \chi_{[0,p^i-1]}(j)\cdot\sumd{k\in\sZ} x_{[j]_{p^{i-m}} - p^{i-m} + kp^{i+s}}.
\]
If $i < m$, we let
\[
\Big(T^{m,s}_i x\Big)_j := 0.
\]
We note that 
\[
\Big(T^{0,0}_i x\Big)_j := \chi_{[0,p^i-1]}(j)\cdot\sumd{k\in\sZ} x_{j + kp^i}.
\]

\begin{Lemma}
\label{LemKM4_7}
\Absit
\begin{itemize}
\item[(i)] Given $a,b,c,d,s,t\geq 0$ such that $b-a\leq d-c\leq b+s \leq d$. Then
\[
T^{a,s}_b \circ T^{c,t}_d = p^{d-b-s} T^{a,d-b+t}_b.
\]
\item[(ii)] Given $a,b,d,s\geq 0$ such that $b+s\leq d$. Then
\[
T^{a,s}_b \circ T^{0,0}_d = T^{a,s}_b.
\]
\end{itemize}

\rm
Ad (i). We may assume $0\leq b-a$. Given $j\in\Z$, we obtain 
\[
\begin{array}{rcl}
\left(T^{a,s}_b T^{c,t}_d x\right)_j
& = & \left(T^{a,s}_b \left(\chi_{[0,p^d-1]}(j')\sumd{k\in\sZ} x_{[j']_{p^{d-c}} - p^{d-c} + kp^{d+t}} \right)_{j'\in\sZ}\right)_j \\
& = & \chi_{[0,p^b-1]}(j)\cdot\sumd{k'\in\sZ} \chi_{[0,p^d-1]}([j]_{p^{b-a}} - p^{b-a} + k'p^{b+s})\cdot \\
&   &                                      \cdot\sumd{k\in\sZ} x_{[[j]_{p^{b-a}} - p^{b-a} + k'p^{b+s}]_{p^{d-c}} - p^{d-c} + kp^{d+t}} \\
& = & \chi_{[0,p^b-1]}(j) p^{d-b-s}\sumd{k\in\sZ} x_{[j]_{p^{b-a}} - p^{b-a} + kp^{d+t}} \\
& = & \left(p^{d-b-s} T^{a,d-b+t}_b x\right)_j.
\end{array}
\]

Ad (ii). The operator $T^{a,s}_b$ is invariant under sequence shifts by $p^{b+s}$, hence under sequence shifts by $p^d$. 
\end{Lemma}

\begin{Lemma} 
\label{LemKM4_7_5}
Given $m\geq 1$, $l\in [1,m]$ and $a\geq 0$, we obtain
\[
\sumd{i\in [0,l-1]} p^{l-1-i} T^{a,l-1-i}_{m-l}\circ (T^{0,0}_{m-i} - T^{1,0}_{m-i}) = T^{a,0}_{m-l} - p^l T^{a,l}_{m-l}. 
\]

\rm
In fact,
\[
\begin{array}{rcl}
\sumd{i\in [0,l-1]} p^{l-1-i} T^{a,l-1-i}_{m-l}\circ (T^{0,0}_{m-i} - T^{1,0}_{m-i})
& \aufgl{(\ref{LemKM4_7} i, ii)} & \sumd{i\in [0,l-1]} p^{l-1-i} T^{a,l-1-i}_{m-l} - \sumd{i\in [0,l-1]} p^{l-i} T^{a,l-i}_{m-l} \\
& = & \sumd{i\in [1,l]} p^{l-i} T^{a,l-i}_{m-l} - \sumd{i\in [0,l-1]} p^{l-i} T^{a,l-i}_{m-l} \\
& = & T^{a,0}_{m-l} - p^l T^{a,l}_{m-l}. \\
\end{array}
\]
\end{Lemma}

\begin{Lemma}
\label{LemKH4_8}
For $x = (x_j)_{j\in\sZ}\in\lfk^\mb{\scr\rm f}(\Z)$ and $l\geq 0$, we denote 
\[
x\ast\zeta_{p^l} := \sum_{j\in\sZ} x_j \zeta_{p^l}^j
\]
and obtain
\[
x\ast\zeta_{p^l} = ((T^{0,0}_l - T^{1,0}_l)x)\ast\zeta_{p^l}.
\]

\rm
We claim that $(T^{1,0}_l x)\ast\zeta_{p^l} = 0$. We may suppose $l\geq 1$ to calculate 
\[
\begin{array}{rcl}
(T^{1,0}_l x)\ast\zeta_{p^l} 
& = & \sum_{j\in\sZ}\chi_{[0,p^l-1]}(j)\sum_{k\in\sZ} x_{[j]_{p^{l-1}} - p^{l-1} + kp^l} \zeta_{p^l}^j \\
& = & \sum_{k\in\sZ} \sum_{h\in [0,p^{l - 1} - 1]}\sum_{i\in [0,p-1]} x_{[h + ip^{l-1}]_{p^{l-1}} - p^{l-1} + kp^l} \zeta_{p^l}^{h+ip^{l-1}} \\
& = & \sum_{k\in\sZ} \sum_{h\in [0,p^{l - 1} - 1]} x_{h - p^{l-1} + kp^l} \zeta_{p^l}^h \sum_{i\in [0,p-1]} \zeta_{p^l}^{ip^{l-1}} \\
& = & 0.\\
\end{array}
\]
\end{Lemma}

\begin{Lemma}
\label{LemKH4_9}
For $l\geq 1$, $s\geq 0$ and $j\in\Z\ohne [0,\phi(p^l) - 1]$, we have $((T^{0,s}_l - T^{1,s}_l)x)_j = 0$.

\rm
In fact,
\[
\begin{array}{rcl}
(T^{1,s}_l x)_{(p-1)p^{l-1} + j}
& = & \sum_{k\in\sZ} x_{[(p-1)p^{l-1} + j]_{p^{l-1}} - p^{l-1} + kp^{l+s}} \\
& = & \sum_{k\in\sZ} x_{j - p^{l-1} + kp^{l+s}} \\
& = & (T^{0,s}_l x)_{(p-1)p^{l-1} + j} \\
\end{array}
\]
for $j\in [0,p^{l-1}-1]$.
\end{Lemma}

\begin{Theorem}
\label{ThKM5}
The image of the absolute Wedderburn embedding is given by
{\rm
\[
\begin{array}{rcl}
(\Z C_{p^n})\omega_{\sZ,p^n}\!\! 
& =   & \left\{\rule[7mm]{0mm}{0mm}\right.\left(\sumd{j\in [0,\phi(p^i) - 1]} x_{i,j}\zeta_{p^i}^j\right)_{\!\! i\in [0,n]},\; x_{i,j}\in\Z\;\left|\rule[7mm]{0mm}{0mm}\right. \\ 
&     & \mb{for $l\in [1,n]$ and $j\in [0,\phi(p^{n-l})-1]$ we have }\; x_{n-l,j} \con_{p^l} \sumd{i\in [0,l-1]} p^{l-1-i} \cdot\\
&     & \cdot \sumd{k\in [1,p-1]} \left(x_{n-i,j - p^{n-l} + kp^{n-1-i}} - (1-\dell_{l,n})x_{n-i,[j]_{p^{n-l-1}} - p^{n-l-1} + kp^{n-1-i}}\right)
        \left.\rule[7mm]{0mm}{0mm}\right\} \\
& =   & \left\{\rule[7mm]{0mm}{0mm}\right.\left(x_i\ast\zeta_{p^i}\right)_{i\in [0,n]},\; x_i\in\lfk^\mb{\scr\rm f}(\Z),\;\ul{x_i}\tm [0,\phi(p^i)-1] \left|\rule[7mm]{0mm}{0mm}\right. \\ 
&     & \mb{ for $l\in [1,n]$ we have }\; x_{n-l} \con_{p^l} \sumd{i\in [0,l-1]} p^{l-1-i}\left(T_{n-l}^{0,l-1-i} - T_{n-l}^{1,l-1-i}\right)x_{n-i} \left.\rule[7mm]{0mm}{0mm}\right\} \\
& \tm & \prodd{i\in [0,n]}\Z[\zeta_{p^i}]. \\
\end{array}
\]}

This system of ties is of triangular shape. In particular, the elementary divisors of $\omega_{\sZ,p^n}$ over $\Z$ are given by $p^i$ with multiplicity $\phi(p^{n-i})$ for $i\in [0,n]$.

\rm
The second description being true if $n = 0$, we perform an induction on $n$. We see by (\ref{PropKM3}) that $\left(x_i\ast\zeta_{p^i}\right)_{i\in [0,n]}\in\prodd{i\in [0,n]}\Z[\zeta_{p^i}]$
is contained in $(\Z C_{p^n})\omega_{\sZ,p^n}$ if and only if the conditions (a) and (b) below hold. To formulate (b), we shall make use of (\ref{LemKH4_8}).
\begin{itemize}
\item[(a)] The tuple $\left(x_i\ast\zeta_{p^i}\right)_{i\in [0,n-1]}$ is contained in $(\Z C_{p^{n-1}})\omega_{\sZ,p^{n-1}}$.

\item[(b)] We have
\[
\left(\left((T^{0,0}_i - T^{1,0}_i)x_n\right)\ast\zeta_{p^i}\right)_{i\in [0,n-1]} - \left(x_i\ast\zeta_{p^i}\right)_{i\in [0,n-1]}\in p(\Z C_{p^{n-1}})\omega_{\sZ,p^{n-1}}.
\]
\end{itemize}

Condition (b) in turn holds if and only if the conditions (ba) and (bb) below hold. To formulate (ba), we use (\ref{LemKH4_9}). To formulate (bb), we shall use that by induction, the description 
is valid in case $n-1$.
\begin{itemize}
\item[(ba)] We have $(T^{0,0}_i - T^{1,0}_i)x_n \con_p x_i$ for $i\in [0,n-1]$.
\item[(bb)] We have
\[
\begin{array}{rl}   
                     & - \left(T^{0,0}_{(n-1)-l} - T^{1,0}_{(n-1)-l}\right)x_n \\
                     & + \sumd{i\in [0,l-1]} p^{l-1-i}\left(T_{(n-1)-l}^{0,l-1-i} - T_{(n-1)-l}^{1,l-1-i}\right)\left(T^{0,0}_{(n-1)-i} - T^{1,0}_{(n-1)-i}\right)x_n \\
\con_{p^{l+1}}       & - x_{(n-1)-l} + \sumd{i\in [0,l-1]} p^{l-1-i}\left(T_{(n-1)-l}^{0,l-1-i} - T_{(n-1)-l}^{1,l-1-i}\right)x_{(n-1)-i} \\
\end{array}
\]
for $l\in [1,n-1]$.
\end{itemize}

By (\ref{LemKM4_7_5}), we may equivalently reformulate to
\begin{itemize}
\item[(bb)] We have
\[
\begin{array}{rl}
0 \con_{p^{l+1}} & - x_{(n-1)-l} +  p^l\left(T_{(n-1)-l}^{0,l} - T_{(n-1)-l}^{1,l}\right)x_n \\
                 & + \sumd{i\in [0,l-1]} p^{l-1-i}\left(T_{(n-1)-l}^{0,l-1-i} - T_{(n-1)-l}^{1,l-1-i}\right)x_{(n-1)-i} \\
\end{array}
\]
for $l\in [1,n-1]$.
\end{itemize}
Shifting the indices $i$ and $l$ by one, we may in turn equivalently reformulate this assertion to 
\begin{itemize}
\item[(bb)] We have
\[
x_{n-l} \con_{p^l} \sumd{i\in [0,l-1]} p^{l-1-i}\left(T_{n-l}^{0,l-1-i} - T_{n-l}^{1,l-1-i}\right)x_{n-i} 
\]
for $l\in [2,n]$.
\end{itemize}

Adjoining condition (ba) in case $i = n-1$, we obtain the necessity of the claimed set of ties.

Conversely, we may use (\ref{LemKM4_7_5}) to see that
\[
(T^{0,0}_{i-1} - T^{1,0}_{i-1})\circ (T^{0,0}_i - T^{1,0}_i) \con_p (T^{0,0}_{i-1} - T^{1,0}_{i-1}) 
\]
for $i\geq 1$, whence the reduction of (bb) to the modulus $p$ suffices to conclude that (ba) holds in all its cases. Moreover, condition (a) follows by induction assumption and by reading
(bb) modulo $p^{l-1}$ for $l\in [2,n]$, thus discarding the summand for $i = 0$. Thus our claimed set of ties is also sufficient.

To see that the second description agrees with the first one, it remains to show that
\[
\left(\left(T_{n-l}^{0,l-1-i} - T_{n-l}^{1,l-1-i}\right)x_{n-i}\right)_j = 0
\]
for $l\in [1,n-1]$, $i\in [0,l-1]$ and $j\in \Z\ohne [0,\phi(p^{n-l})-1]$. But this follows from (\ref{LemKH4_9}).
\end{Theorem}

\begin{quote}
\begin{footnotesize}
\begin{Remark}
\label{RemKM4_7_7}\rm
It might be worthwhile to try to give an ad hoc proof of (\ref{ThKM5}) by verification of the ties on the image of the canonical group basis of $\Z C_{p^n}$ and by a comparison of indices -- which 
we have not attempted to do. We preferred to proceed in the straightforward manner as above since, in this way, the role of the Kervaire-Murthy pullback (\ref{PropKM3}) remains visible. 
\end{Remark}

\begin{Example}
\label{ExKM5_5}
\rm
Let $p = 3$. We obtain 
\[
\begin{array}{rcl}
\Z C_3    & \lraisoa{\omega_{\ssZ,3}}  & \Big\{ \left(\sum_{k\in [0,1]} x_{1,k}\zeta_3^k\right) \ti x_{0,0} \; \Big| \; x_{0,0} \con_3 x_{1,0} + x_{1,1} \Big\} \; \tm \; \Z[\zeta_3]\ti\Z \\
\Z C_9    & \lraisoa{\omega_{\ssZ,9}}  & \Big\{\left(\sum_{k\in [0,5]} x_{2,k}\zeta_9^k\right) \ti \left(\sum_{k\in [0,1]} x_{1,k}\zeta_3^k\right) \ti x_{0,0} \;\Big| \\
& & x_{1,0} \con_3 x_{2,0} + x_{2,3} - x_{2,2} - x_{2,5}  \\
& & x_{1,1} \con_3 x_{2,1} + x_{2,4} - x_{2,2} - x_{2,5}  \\
& & x_{0,0} \con_9 3(x_{2,2} + x_{2,5}) + (x_{1,0} + x_{1,1}) \Big\} \; \tm \; \Z[\zeta_9]\ti\Z[\zeta_3]\ti\Z \\
\Z C_{27} & \lraisoa{\omega_{\ssZ,27}} & \Big\{\left(\sum_{k\in [0,17]} x_{3,k}\zeta_{27}^k\right)\ti\left(\sum_{k\in [0,5]} x_{2,k}\zeta_9^k\right)
                                               \ti\left(\sum_{k\in [0,1]} x_{1,k}\zeta_3^k\right) \ti x_{0,0} \;\Big| \\
& & x_{2,0} \con_3 x_{3,0} + x_{3,9}  - x_{3,6} - x_{3,15}  \\
& & x_{2,1} \con_3 x_{3,1} + x_{3,10} - x_{3,7} - x_{3,16}  \\
& & x_{2,2} \con_3 x_{3,2} + x_{3,11} - x_{3,8} - x_{3,17}  \\
& & x_{2,3} \con_3 x_{3,3} + x_{3,12} - x_{3,6} - x_{3,15}  \\
& & x_{2,4} \con_3 x_{3,4} + x_{3,13} - x_{3,7} - x_{3,16}  \\
& & x_{2,5} \con_3 x_{3,5} + x_{3,14} - x_{3,8} - x_{3,17}  \\
& & x_{1,0} \con_9 3(x_{3,6} + x_{3,15} - x_{3,8} - x_{3,17}) + (x_{2,0} + x_{2,3} - x_{2,2} - x_{2,5}) \\
& & x_{1,1} \con_9 3(x_{3,7} + x_{3,16} - x_{3,8} - x_{3,17}) + (x_{2,1} + x_{2,4} - x_{2,2} - x_{2,5}) \\
& & x_{0,0} \con_{27} 9(x_{3,8} + x_{3,17}) + 3(x_{2,2} + x_{2,5}) + (x_{1,0} + x_{1,1})  \Big\} \\
& & \tm\; \Z[\zeta_{27}]\ti\Z[\zeta_9]\ti\Z[\zeta_3]\ti\Z. \\
\end{array}
\]
\end{Example}
\end{footnotesize}
\end{quote}

\subsection{Kervaire-Murthy ties for $\Z C_m$}
\label{SubSecKerMurTiesm}

Let $m\geq 1$. We maintain the notation of section \ref{SubSecInvm}.
Suppose given an inclusion of commutative $\Z$-orders $\Lambda\tm\Gamma$ that has, as an inclusion of abelian groups, a cokernel $\Gamma/\Lambda$ which is finite 
as a set. Let the {\it naive localization} $\Lambda_{[p],\Gamma}$ be the kernel of the composition $\Gamma\lra\Gamma/\Lambda\lra(\Gamma/\Lambda)_{(p)}$ 
(cf.\ [K\"u 99, D.2]). Since $\Lambda_{[p],\Gamma} = \Gamma\cap\Lambda_{(p)}$, intersected as subsets of $\Gamma_{(p)}$, the naive localization $\Lambda_{[p],\Gamma}$ is a suborder of 
$\Gamma$ of index $(\#(\Gamma/\Lambda))[p]$. We have $\Lambda_{(p)} = (\Lambda_{[p],\Gamma})_{(p)}$. Moreover, 
\[
\Lambda = \Cap_p\Lambda_{[p],\Gamma}.
\]
Note that if $\Lambda_{(p)} = \Gamma_{(p)}$, then $\Lambda_{[p],\Gamma} = \Gamma$.

\begin{Lemma}
\label{LemCP9_1}
Given inclusions of $\Z$-orders $\Lambda\tm\Gamma$ and $\Lambda'\tm\Gamma'$, we obtain an equality
\[
(\Lambda\ts\Lambda')_{[p],\Gamma\ts\Gamma'} = \Lambda_{[p],\Gamma}\ts\Lambda'_{[p],\Gamma'}.
\]
as subsets of $\Gamma\ts\Gamma'$.

\rm
In fact, flatness yields
\[
\begin{array}{rcl}
\Lambda_{[p],\Gamma}\ts \Lambda'_{[p],\Gamma'}
& = & (\Lambda_{[p],\Gamma}\ts\Lambda'_{(p)})\cap (\Lambda_{[p],\Gamma}\ts\Gamma') \\
& = & (\Lambda\ts\Lambda')_{(p)}\cap (\Lambda_{(p)}\ts\Gamma') \cap (\Gamma\ts\Gamma') \\
& = & (\Lambda\ts\Lambda')_{[p],\Gamma\ts\Gamma'}.\\
\end{array}
\]
\end{Lemma}

Letting 
\[
\begin{array}{rcrcc}
(\Lambda\tm\Gamma)   & := & \Big((\Z C_{m[p]})\omega_{\sZ,m[p]}   & \hra & \prodd{e|m[p]}\Z[\zeta_e]\Big) \\ 
(\Lambda'\tm\Gamma') & := & \Big((\Z C_{m[p']})\omega_{\sZ,m[p']} & \hra & \prodd{f|m[p']}\Z[\zeta_f]\Big), \\
\end{array}
\]
we obtain
\[
\begin{array}{rcl}
(\Lambda\ts\Lambda')_{[p],\Gamma\ts\Gamma'} 
&\aufgl{(\ref{LemCP9_1})} & \Lambda_{[p],\Gamma}\ts \Lambda'_{[p],\Gamma'} \\
&\aufgl{(\ref{CorCP9_0_5})} & \Lambda\ts\Gamma'\\
& = & \prodd{f|m[p']} \Lambda\ts \Z[\zeta_f] \\
& \tm & \prodd{f|m[p']}\left(\prodd{e|m[p]}\Z[\zeta_e]\right)\ts\Z[\zeta_f] \\
& = & \prodd{f|m[p']}\;\prodd{e|m[p]}\Z[\zeta_{ef}]. \\
\end{array}
\]

\begin{quote}
\begin{footnotesize}
We could not do better than to argue with a $\Z$-linear basis of $\Z[\zeta_f]$. This lead to the `non-canonical' representation of an element of $\Z[\zeta_d]$, $d|m$, that is used in the following
\end{footnotesize}
\end{quote}

\begin{Proposition}
\label{PropCP9_2}
Suppose given an element 
\[
\left(\sumd{(j_p)_{p|d}\in\prod_{p|d} [0,\phi(d[p]) - 1]} a_{d,(j_p)_{p|d}} \zeta_d^{\sum_{p|d} s_{d,p} j_p d[p']}\right)_{\!\!\! d|m}\in \prodd{d|m}\Z[\zeta_d],
\] 
where $a_{d,(j_p)_{p|d}}\in\Z$. This element is contained in the image of $\omega_{\sZ,m}$ if and only if for each $p|m$, for each $f | m[p']$ and for each tuple $(j_q)_{q|f}$, 
$j_q\in [0,\phi(d[q])-1]$, we have
\[
\left(\sumd{j_p\in [0,\phi(e)-1]} a_{ef,\; j_p\ti (j_q)_{q|f}}\zeta_e^{j_p}\right)_{\!\! e|m[p]}\in (\Z C_{m[p]})\omega_{\sZ,m[p]},
\]
where $j_p\ti (j_q)_{q|f}$ denotes the tuple that has entry $j_p$ at position $p$ and $j_q$ at position $q$ for $q|f$. Thus we may employ (\ref{ThKM5}).
\end{Proposition}

\begin{quote}
\begin{footnotesize}
\begin{Example}
\label{ExCP9_3}
\rm 
We have $s_{6,2} = 1$, $s_{6,3} = 2$ and $s_{12,2} = 3$, $s_{12,3} = 1$. The element
\[
\begin{array}{l}
(a_{1,(\bullet,\bullet)})\ti (a_{2,(0,\bullet)})\ti (a_{3,(\bullet,0)} + a_{3,(\bullet,1)}\zeta_3) \ti (a_{4,(0,\bullet)} + a_{4,(1,\bullet)}\zeta_4)\ti (a_{6,(0,0)} + a_{6,(0,1)}\zeta_6^4) \\
\ti (a_{12,(0,0)} + a_{12,(0,1)}\zeta_{12}^4 + a_{12,(1,0)}\zeta_{12}^9 + a_{12,(1,1)}\zeta_{12}), \\
\end{array}
\]
the symbol $\bullet$ indicating a non-existing entry, is contained in $(\Z C_{12})\omega_{\ssZ,12}$ if and only if
\[
\begin{array}{lclclcl}
(a_{1,(\bullet,\bullet)}) & \ti & (a_{2,(0,\bullet)})       & \ti & (a_{4,(0,\bullet)} + a_{4,(1,\bullet)}\zeta_4) & \in & (\Z C_4)\omega_{\ssZ,4} \\
(a_{3,(\bullet,0)})       & \ti & (a_{6,(0,0)})             & \ti & (a_{12,(0,0)} + a_{12,(1,0)}\zeta_4)           & \in & (\Z C_4)\omega_{\ssZ,4} \\
(a_{3,(\bullet,1)})       & \ti & (a_{6,(0,1)})             & \ti & (a_{12,(0,1)} + a_{12,(1,1)}\zeta_4)           & \in & (\Z C_4)\omega_{\ssZ,4} \\
                          &     & (a_{1,(\bullet,\bullet)}) & \ti & (a_{3,(\bullet,0)} + a_{3,(\bullet,1)}\zeta_3) & \in & (\Z C_3)\omega_{\ssZ,3} \\
                          &     & (a_{2,(0,\bullet)})       & \ti & (a_{6,(0,0)} + a_{6,(0,1)}\zeta_3)             & \in & (\Z C_3)\omega_{\ssZ,3} \\
                          &     & (a_{4,(0,\bullet)})       & \ti & (a_{12,(0,0)} + a_{12,(0,1)}\zeta_3)           & \in & (\Z C_3)\omega_{\ssZ,3} \\
                          &     & (a_{4,(1,\bullet)})       & \ti & (a_{12,(1,0)} + a_{12,(1,1)}\zeta_3)           & \in & (\Z C_3)\omega_{\ssZ,3}. \\
\end{array}
\]
By (\ref{ThKM5}), an element $(b_{1,0}) \ti (b_{2,0}) \ti (b_{4,0} + b_{4,1}\zeta_4)$ is contained in $(\Z C_4)\omega_{\ssZ,4}$ if and only if
\[
\begin{array}{rcl}
b_{2,0} & \con_2 & b_{4,0} - b_{4,1} \\
b_{1,0} & \con_4 & 2 b_{4,1} + b_{2,0}, \\
\end{array}
\]
and an element $(b_{1,0}) \ti (b_{3,0} + b_{3,1}\zeta_3)$ is contained in $(\Z C_3)\omega_{\ssZ,3}$ if and only if
\[
\begin{array}{rcl}
b_{1,0} \con_3 b_{3,0} + b_{3,1}. \\
\end{array}
\]
\end{Example}
\end{footnotesize}
\end{quote}

\begin{Remark}
\label{RemCP11}
\rm
{\sc Kleinert} gives a system of ties that describes the image of the absolute Wedderburn embedding $\Z C_m\lraa{\omega_{\ssZ,m}}\prod_{d|m} \Z[\zeta_d]$ in terms of certain 
prime ideals of $\Z[\zeta_d]$, $d|m$, in case $m$ is squarefree [Kl 81, p.\ 550].
\end{Remark}
\section{References}

\begin{footnotesize}
{\sc G.\ E.\ Andrews}\upl
\begin{itemize}
\item[{[A 76]}] {\it The Theory of Partitions,} Encyclopedia of Math.\ and Appl., Vol.\ 2, Addison-Wesley, 1976.
\end{itemize}

{\sc Buenos Aires Cyclic Homology Group (J.\ A.\ Guccione, J.\ J.\ Guccione, M.\ J.\ Redondo, A.\ Solotar, O.\ E.\ Villamayor)}\upl
\begin{itemize}
\item[{[BACH 91]}] {\it Cyclic Homology of Algebras with One Generator,} $K$-Theory 5, p.\ 51-69, 1991. 
\end{itemize}

{\sc M.A.\ Kervaire, M.P.\ Murthy}\upl
\begin{itemize}
\item[{[KM 77]}] {\it On the projective class group of cyclic groups of prime power order,} Comm.\ Math.\ Helv.\ 52, p.\ 415-452, 1977. 
\end{itemize}

{\sc M.\ K\"unzer}\upl
\begin{itemize}
\item[{[K\"u 99]}] {\it Ties for the integral group ring of the symmetric group,} thesis, \\ 
                   http://www.mathematik.uni-bielefeld.de/$\scm\sim$kuenzer, Bielefeld, 1999.
\end{itemize}

{\sc Kleinert, E.}\upl
\begin{itemize}
\item[{[Kl 81]}] {\it Einheiten in $\Z[D_{2m}]$,} J.\ Num.\ Th.\ 13, p.\ 541-561, 1981. 
\end{itemize}

{\sc M.\ Larsen, A.\ Lindenstrauss}\upl
\begin{itemize}
\item[{[LL 92]}] {\it Cyclic homology of Dedekind domains,} $K$-Theory 6, p.\ 301-334, 1992. 
\end{itemize}

{\sc J.\ Neukirch}\upl
\begin{itemize}
\item[{[N 91]}] {\it Algebraische Zahlentheorie,} Springer, 1991.
\end{itemize}

{\sc Plesken, W.}\upl
\begin{itemize}
\item[{[P 80]}] {\it Gruppenringe \"uber lokalen Dedekindbereichen}, Habilitationsschrift, Aachen, 1980.
\end{itemize}

\vspace*{1cm}

\begin{flushleft}
Harald Weber\\
Mathematisches Institut B\\
3. Lehrstuhl\\
Pfaffenwaldring 57\\ 
70569 Stuttgart\\
harald@poolb.mathematik.uni-stuttgart.de\\
\ \\
Matthias K\"unzer\\
Fakult\"at f\"ur Mathematik\\
Universit\"at Bielefeld\\
Postfach 100131\\
D-33501 Bielefeld\\
kuenzer@mathematik.uni-bielefeld.de\\
\end{flushleft}
\end{footnotesize}
\end{document}